\documentclass[12pt]{article}
\usepackage{subfigure}
\usepackage{tikz}
\usepackage{pstricks, amssymb, xr, multind, multicol,textcomp,amsthm}
\def\hfl#1{\smash{\mathop{\hbox to 10mm{\rightarrowfill}}\limits^{\textstyle
#1}}}
\newtheorem{proposition}[equation]{Proposition} 
\newtheorem{corollary}[equation]{Corollary} 
\newtheorem{theorem}[equation]{Theorem} 
\newtheorem{exa}[equation]{Example} 
\newtheorem{ex}[equation]{Exercise} 
\newtheorem{s-ex}[equation]{Side-exercise} 
\newtheorem{exas}[equation]{Examples} 
\newtheorem{lemma}[equation]{Lemma} 
\newtheorem{sublemma}[equation]{Sublemma} 
\newtheorem{remar}[equation]{Remark}
\newtheorem*{remarnn}{Remark}
\newtheorem{remars}[equation]{Remarks}
\newtheorem{nota}[equation]{Notation}
\newtheorem{sremar}[equation]{Side-remark} 
\newtheorem{definitio}[equation]{Definition}

\newenvironment{notation}{\begin{nota} \rm }{\end{nota}} 
\newenvironment{remark}{\begin{remar} \rm }{\end{remar}}

\newenvironment{definition}{\begin{definitio} \rm }{\end{definitio}}

\newcommand{\HH}{{\mathbb{H}}} 
\newcommand{\KK}{\mathbb{K}}

\newcommand{\CS}{{\cal S}}

\newcommand{\ZZ}{\mathbb{Z}}

\newcommand{\RR}{\mathbb{R}} 
\newcommand{\QQ}{\mathbb{Q}} 
\newcommand{\CC}{\mathbb{C}}

\newcommand{\bp}{\noindent {\sc Proof: }} 
\newcommand{\eop}{\nopagebreak \hspace*{\fill}{$\diamond$} \medskip} 
\newcommand{\neghopf}{\begin{tikzpicture}[x=-30,scale=.5] \useasboundingbox (-1,-.5) rectangle (.7,.5);
\draw[->] (0,0) -- (0,.2) .. controls (0,.4) and (-.15,.5).. (-.3,.5) .. controls(-.45,.5) and (-.6,.2) .. (-.6,0);
\draw[->] (-.6,0) .. controls (-.6,-.2) and (-.45,-.5) .. (-.3,-.5) .. controls (-.15,-.5) and (0,-.35) .. (0,-.2);
\draw[->] (-.1,.2) .. controls (-.2,.2) and (-.25,.15) .. (-.25,.05) .. controls (-.25,-.05) and (-.2,-.1) .. (0,-.1) .. controls (.2,-.1) and (.25,-.05) .. (.25,.05) .. controls (.25,.15) and (.2,.2) .. (.1,.2);
\end{tikzpicture}} 
\newcommand{\poshopf}{\begin{tikzpicture}[scale=.5] \useasboundingbox (-1,-.5) rectangle (1,.5);
\draw[->] (0,0) -- (0,.2) .. controls (0,.4) and (-.15,.5).. (-.3,.5) .. controls(-.45,.5) and (-.6,.2) .. (-.6,0);
\draw[->] (-.6,0) .. controls (-.6,-.2) and (-.45,-.5) .. (-.3,-.5) .. controls (-.15,-.5) and (0,-.35) .. (0,-.2);
\draw[->] (-.1,.2) .. controls (-.2,.2) and (-.25,.15) .. (-.25,.05) .. controls (-.25,-.05) and (-.2,-.1) .. (0,-.1) .. controls (.2,-.1) and (.25,-.05) .. (.25,.05) .. controls (.25,.15) and (.2,.2) .. (.1,.2);
\end{tikzpicture}}
\begin{document} 
\title{On homotopy invariants of combings of three-manifolds}
\author{Christine Lescop \thanks{Institut Fourier, UJF Grenoble, CNRS}}
\maketitle
\begin{abstract} 
Combings of oriented compact $3$-manifolds are homotopy classes of nowhere zero vector fields in these manifolds.
A first known invariant of a combing is its Euler class, that is the Euler class of the normal bundle to a combing representative in the tangent bundle of the $3$-manifold $M$. It only depends on the Spin$^c$-structure represented by the combing.
When this Euler class is a torsion element of $H^2(M;\ZZ)$, we say that the combing is a torsion combing. Gompf introduced a $\QQ$-valued invariant $\theta_G$ of torsion combings of closed $3$-manifolds that distinguishes all combings that represent a given Spin$^c$-structure. This invariant provides a grading of the Heegaard Floer homology $\widehat{HF}$ for manifolds equipped with torsion Spin$^c$-structures.
We give an alternative definition of the Gompf invariant and we express its variation as a linking number. We also define a similar invariant $p_1$ for combings of manifolds bounded by $S^2$. We show that the $\Theta$-invariant, that is the simplest configuration space integral invariant of rational homology spheres, is naturally an invariant of combings of rational homology balls, that reads $(\frac14p_1 + 6 \lambda)$ where $\lambda$ is the Casson-Walker invariant.
The article also includes a mostly self-contained presentation of combings.
\end{abstract}

\maketitle
\section{Introduction}

\subsection{General introduction}
In this article, $M$ is an oriented connected compact smooth $3$-manifold.
The boundary $\partial M$ of $M$ is either empty or identified with the unit sphere $S^2$ of $\RR^3$.
In this latter case, a neighborhood $N(\partial M)$ of $\partial M$ in $M$ is identified with a neighborhood of $S^2$ in the unit ball of $\RR^3$.
 The tangent bundle of $M$ is denoted by $TM$, and the unit tangent bundle of $M$ is denoted by $UM$. Its fiber is $U_mM =(T_mM \setminus \{0\})/\RR^{\ast +}$. 
All parallelizations of $M$ are assumed to coincide with the parallelization induced by the standard parallelization $\tau_s$ of $\RR^3$ over $N(\partial M)$, and all sections of $UM$ are assumed to be constant with respect to this parallelization over $N(\partial M)$. Homotopies of parallelizations or sections satisfy these assumptions at any time. When $\partial M=\emptyset$, the parallelizations of $M$ also induce the orientation of $M$.

A {\em combing\/} of $M$ is a homotopy class of such sections of $UM$.
According to Turaev \cite{Tu}, a {\em Spin$^c$-structure\/} on $M$ may be seen as an equivalence class of sections of $UM$,
where two sections are in the same class if and only if they are homotopic over the complement of a point that sits in the interior of $M$.

For $\KK=\ZZ$ or $\QQ$, a {\em $\KK$-sphere or (integral or rational) homology sphere\/} (resp. a {\em $\KK$-ball\/}) is a smooth, compact, oriented $3$-manifold with the same $\KK$-homology as the sphere $S^3$ (resp. as a point).

In this mostly self-contained article, we study the combings of $M$, that are homotopy classes of sections of $UM$. We describe their classification, and some of their invariants. We first describe the first known homotopy invariant of a combing, that is the Euler class, in terms of links. The {\em Euler class\/} of a combing is the Euler class of the normal bundle to a combing representative in $TM$. It only depends on the Spin$^c$-structure induced by the combing.
When this Euler class is a torsion element of $H^2(M,\partial M;\ZZ)$, we say that the combing is a {\em torsion combing\/}.
We introduce a rational invariant $p_1$ of torsion combings of $M$.
When $M$ is \emph{closed} (i.e. compact, without boundary), we show that the invariant $p_1$ coincides with an invariant $\theta_G$ defined by Gompf in \cite{Go}.
For a combing that extends to a parallelization, the invariant $p_1$ coincides with the Hirzebruch defect (or Pontrjagin number) of the parallelization, studied in \cite{hirzebruchEM, km, lesconst, lesmek}.
In general, we express the variation of $p_1$ in terms of linking numbers.
We also give a homogeneous self-contained definition of an invariant $\Theta$ of combings of rational homology balls from configuration spaces, and we show that this $\Theta$-invariant reads $(\frac14p_1 + 6 \lambda)$ where $\lambda$ is the Casson-Walker invariant normalized as in \cite{akmc,mar} for $\ZZ$-spheres and as $\frac{\lambda_W}{2}$ for $\QQ$-spheres, where $\lambda_W$ is the Walker normalization in \cite{wal}.

\subsection{Conventions and notations}
Unless otherwise mentioned, all manifolds are oriented. Boundaries are oriented by the outward normal first convention.
Products are oriented by the order of the factors. More generally, unless otherwise mentioned, the order of appearance of coordinates or parameters orients chains or manifolds.
The fiber $N_x(A)$ of the normal bundle $N(A)$ to an oriented submanifold $A$ at $x \in A$ is oriented so that $N_x(A)$ followed by the tangent bundle $T_x(A)$ to $A$ at $x$ induces the orientation of the ambient manifold. The orientation of $N_x(A)$ is a \emph{coorientation} of $A$ at $x$. The transverse preimage of a submanifold under a map $f$ is oriented so that $f$ preserves the coorientations.
The transverse intersection of two submanifolds $A$ and $B$ in a manifold $M$ is oriented so that the normal bundle to $A\cap B$ is $(N(A) \oplus N(B))$, fiberwise.
If the two manifolds are of complementary dimensions, then the sign of an intersection point is $+1$ if the orientation of its normal bundle coincides with the orientation of the ambient space, that is if $T_xM=N_xA \oplus N_xB$ (as oriented vector spaces), this is equivalent to $T_xM=T_xA \oplus T_xB$. Otherwise, the sign is $-1$. If $A$ and $B$ are compact and if $A$ and $B$ are of complementary dimensions in $M$, their {\em algebraic intersection\/} is the sum of the signs of the intersection points, it is denoted by
$\langle A, B \rangle_M$.
The {\em linking number\/} of two rationally null-homologous disjoint links in a $3$-manifold is the algebraic intersection of a rational chain bounded by one of the links and the other one.

\subsection{Expanded introduction}
Let us now be more explicit in order to state the main results precisely. The assertions below will be justified in Subsections~\ref{subPontbiscons} and \ref{subdetint}.
Recall that any smooth compact oriented $3$-manifold $M$
can be equipped with a parallelization $\tau\colon M \times \RR^3 \rightarrow TM$.
When such a parallelization $\tau$ of $M$ is given, two sections $X$ and $Y$ of $UM$ induce a map $(X,Y)\colon M \rightarrow S^2 \times S^2$. Two sections $X$ and $Y$ are said to be {\em transverse\/} if the induced maps $(X,Y)$ and $(X,-Y)$ are transverse to the diagonal of $S^2 \times S^2$, that is if their images are. This is generic and independent of $\tau$.
For two transverse sections $X$ and $Y$, let $L_{X=Y}$ be the preimage of the diagonal of $S^2$ under the map $(X,Y)$. Thus $L_{X=Y}$ is an oriented link in the interior of $M$. It is cooriented by the fiber of the normal bundle to the diagonal of $(S^2)^2$.

The Spin$^c$-structures of $M$ form an affine space $\CS(M)$ with translation group $H^2(M,\partial M;\pi_2(S^2))$.
The Poincar\'e duality isomorphism $P \colon H^2(M,\partial M;\ZZ) \rightarrow H_1(M;\ZZ)$ identifies the translation group of $\CS(M)$ with $H_1(M;\ZZ)$. Let $[X]^c$ denote the Spin$^c$-structure of $M$ represented by a section $X$ of $UM$.
Then, for any two transverse sections $X$ and $Y$ of $UM$, the difference $([X]^c-[Y]^c) \in H_1(M;\ZZ)$ is the class of $L_{X=-Y}$ in $H_1(M;\ZZ)$.

The Euler class of a combing $[X]$ represented by a section $X$, is the Euler class of the normal bundle $TM/\RR X$. It is denoted as $e(X^{\perp})$, it belongs to $H^2(M;\ZZ)$ (here, $H^2(M;\ZZ)=H^2(M,\partial M;\ZZ)$)  and satisfies
$$P(e(X^{\perp}))=[X]^c-[-X]^c$$
so that for two combings $[X]$ and $[Y]$,  $$P(e(X^{\perp})-e(Y^{\perp}))=2([X]^c-[Y]^c).$$

A {\em torsion combing\/} of $M$ is a combing whose Euler class is a torsion element of $H^2(M,\partial M;\ZZ)$. A {\em torsion section\/} of $UM$ is a section that represents a torsion combing.

There is a natural transitive action of $\pi_3(S^2)=\ZZ$ on the combings of $M$ that belong to a given Spin$^c$-structure. This action is free for {\em torsion Spin$^c$-structures\/}, that are Spin$^c$-structures represented by torsion combings. In general, the action of $\pi_3(S^2)$ induces a free transitive action of $\ZZ /\langle e(X^{\perp}),H_2(M,\partial M;\ZZ)\rangle$ on the combings of $M$ that belong to a given Spin$^c$-structure $[X]^c$.

We prove the following theorem in Subsection~\ref{subpfcartorsec}, with elementary arguments.
\begin{theorem}
\label{thmcartorsec}
 Let $X$ be a fixed section of $UM$.
Two sections $Y$ and $Y^{\prime}$ of $UM$ transverse to $X$ represent the same Spin$^c$-structure if and only if the links $L_{Y=-X}$ and $L_{Y^{\prime}=-X}$ are homologous.

If $X$ is a torsion section,
then a section $Y$ of $UM$ transverse to $X$ is a torsion section if and only if the links $L_{Y=-X}$ and $L_{Y=X}$ are rationally null-homologous in $M$.

If $X$ is a torsion section, then two torsion sections $Y$ and $Y^{\prime}$ of $UM$ transverse to $X$ represent the same combing if and only if the links $L_{Y=-X}$ and $L_{Y^{\prime}=-X}$ are homologous, and $lk(L_{Y=-X},L_{Y=X})=lk(L_{Y^{\prime}=-X},L_{Y^{\prime}=X})$.
\end{theorem}
This theorem is a variant of a Pontrjagin theorem recalled in Subsection~\ref{subPontbiscons} that treats the case when $X$ extends to a parallelization.
It might be already known. I thank Patrick Massot for pointing out to me that Dufraine proved similar results in \cite{dufraine}.

The first Pontrjagin class induces a canonical map $p_1$ from the set of parallelization homotopy classes of $M$ to $\ZZ$. When $\partial M=\emptyset$, the map $p_1$, denoted as $\delta(M,.)$, is studied by Hirzebruch in \cite[\S 3.1]{hirzebruchEM}, and Kirby and Melvin study $p_1$ under the name {\em Hirzebruch defect\/} in \cite{km}, and they denote it as $h$, there. This map $p_1$ is studied in \cite{lesconst,lesmek} when $M$ is a $\QQ$-ball. The definition of $p_1$ and some of its properties are recalled in Subsection~\ref{subrecpone}.

The main original result of this article is the following theorem that is proved in Subsection~\ref{subpfmth}.

\begin{theorem}
\label{thmponegen}
There exists a unique map $$p_1\colon \{\mbox{Torsion combings of } M  \} \rightarrow \QQ$$ such that
\begin{itemize}
\item if the combing $[X]$ extends as a parallelization $\tau$, then $p_1([X])=p_1(\tau)$, and
\item for any two transverse torsion sections $X$ and $Y$ of $UM$,
$$p_1([Y])-p_1([X])=4lk(L_{X=Y},L_{X=-Y}).$$
\end{itemize}
The map $p_1$ satisfies the following properties:
\begin{itemize}
\item For any combing $[X]$, $p_1([X])=p_1([-X])$.
\item The restriction of $p_1$ to any torsion Spin$^c$-structure is injective.
\end{itemize}
\end{theorem}

The variation of $p_1$ under simple operations on torsion combings is presented in Subsection~\ref{subonemore}. 

The image of $p_1$ is determined by the following theorem that is proved in Subsection~\ref{subpfmth}.

Let $\ell \colon \mbox{Torsion}(H_1(M;\ZZ)) \rightarrow \QQ/\ZZ$ denote the {\em self-linking number\/} (the linking number of a representative and one of its parallels). View an element $\overline{a}$ of $\QQ/\ZZ$ as its class $(a+\ZZ)$ in $\QQ$ so that $4 \ell(\mbox{Torsion}(H_1(M;\ZZ))$ is a subset of $\QQ$, invariant by translation by $4$.
\begin{theorem}
\label{thmimcom}
Let $\tau$ be a parallelization of $M$ inducing a combing $X$. For any torsion combing $Y$,
$$p_1(Y)\in (p_1(\tau)-4 \ell([L_{Y=-X}])).$$
$$p_1(\{\mbox{\rm Torsion combings}\})=p_1(\tau) -4 \ell(\mbox{\rm Torsion}(H_1(M;\ZZ)).$$
\end{theorem}
Here $p_1(\tau)$ is an integer whose parity is determined in Theorem~\ref{thmim}.
Note that the image of $p_1$ is not an affine space in general.

In Subsection~\ref{secGo}, we  prove that the invariant $p_1$ coincides with the Gompf invariant when $\partial M=\emptyset$.
The  Gompf invariant is denoted by $\theta$ in \cite{Go}, and it is denoted by $\theta_G$ in this article to prevent confusion with $\Theta$.

In \cite[Section 2.6]{OS1}, Ozsv\'ath and Szab\'o associate a Spin$^c$-structure with a generator ${\bf x}$ of the Heegaard Floer homology $\widehat{HF}$. Gripp and Huang refine this process in \cite{GH} in order to associate a combing $\tilde{gr}({\bf x})$ with such a generator ${\bf x}$, and they relate the Gompf invariant to the absolute $\QQ$-grading $\overline{gr}$ of Ozsv\'ath and Szab\'o for the Heegaard Floer homology of $3$-manifolds equipped with torsion Spin$^c$ structures in \cite{OS2}. According to \cite[Corollary 4.3]{GH}, $\overline{gr}({\bf x})=\frac{2+\theta_G(\tilde{gr}({\bf x}))}4$.

The work of Witten \cite{witten} pioneered the introduction of many $\QQ$-sphere invariants, and Witten's insight into the perturbative expansion of Chern-Simons theory led Kontsevich to outline a construction of invariants associated with graph configuration spaces in \cite{ko}. In \cite{kt}, G. Kuperberg and D. Thurston applied the Kontsevich scheme to show the existence of such an invariant $Z_{KKT}$ of $\QQ$-spheres that is equivalent to the LMO invariant of Le, Murakami and Ohtsuki \cite{lmo} for integral homology spheres.
This invariant $Z_{KKT}$ is in fact a graded invariant of parallelized $\QQ$-balls $M$.
Its degree one part is called the $\Theta$-invariant. Let us denote it by $\Theta_{KKT}$. 
For a $\QQ$-ball $M$ equipped with a parallelization $\tau$, the invariant $\Theta_{KKT}(M,\tau)$ is
the sum of $6 \lambda(M)$ and $\frac{p_1(\tau)}{4}$, 
where $\lambda$ is the Casson-Walker invariant,
according to a Kuperberg-Thurston theorem \cite{kt} generalized to rational homology spheres in \cite[Theorem 2.6 and Section 6.5]{lessumgen}.

In Section~\ref{sectheta}, we give a self-contained homogeneous definition of an invariant $\Theta$ of combings $[X]$ in a $\QQ$-sphere $M$ from an algebraic intersection in a two-point configuration space. This invariant satisfies the same variation formula as $\frac14 p_1$ so that $\Theta(M,X)-\frac{p_1([X])}{4}$ only depends on the $\QQ$-sphere $M$.
Furthermore, when $X$ is the first vector of a trivialization $\tau$, it is easy to see that the definition of $\Theta(M,X)$ agrees with the definition of $\Theta_{KKT}(M,\tau)$ as an algebraic intersection of three chains in a two-point configuration space that can be found in \cite[Section 6.5]{lessumgen} and in \cite[Theorem 2.14]{lesbonn}
so that 
$$\Theta(M,X)=6\lambda(M)+\frac14 p_1(X).$$

\section{Combings}
\setcounter{equation}{0}

\subsection{Generalization of a Pontrjagin construction in dimension 3}
\label{subPontbiscons}

\begin{lemma}
\label{lemLXYwd}
Combings are generically transverse. For two transverse sections $X$ and $Y$ of\/ $UM$, the homology classes of\/ $L_{X=Y}$ and $L_{X=-Y}$ only depend on the Spin$^c$-structures $[X]^c$ and $[Y]^c$ represented by $X$ and $Y$.
\end{lemma}
\bp
When $X$ extends as a parallelization, this parallelization identifies $UM$ with $M\times S^2$, then $Y$ may be seen as a map from $M$ to $S^2$, and a homotopy of $Y$ is a map from $[0,1]\times M$ to $S^2$, for which $X$ is a regular value, generically. In particular, the preimage of $X$ under such a homotopy $h$
yields a cobordism from $L_{Y_{0}=X}$ and $L_{Y_{1}=X}$, and the homology class of $L_{Y=X}$ only depends on the homotopy class of $Y$, when $X$ is fixed.
Since any $X$ locally extends as a parallelization, the local transversality arguments hold for any $X$ so that the above proof may be
adapted to any $X$ by using a homotopy $(Y_t,X)$ valued in $S^2 \times S^2$ (with respect to some reference trivialization) and the preimage of the diagonal under this homotopy.
Similarly, the homology class of $L_{Y=-X}$ only depends on the homotopy classes of $X$ and $Y$. Since the homology classes of $L_{Y=-X}$
and $L_{Y=X}$ are unchanged under a modification of $X$ or $Y$ supported in a ball, they only depend on $[X]^c$ and $[Y]^c$.
\eop

Let $X$ be a section of $UM$. Equip $M$ with a Riemannian structure (all of them are homotopic).
These two assumptions hold for the rest of the subsection.

Let $NL$ be the normal bundle to a link $L$ in $M$.
Let $S(NL,(-X)^{\perp})$ denote the space of homotopy classes of sections of the $S^1$--bundle over $L$ whose fiber over $x$ is the space of orientation-preserving linear isometries from the fiber $N_xL=(T_xL)^{\perp} \cong T_xM/T_xL$ of $NL$ to $(-X(x))^{\perp}\cong T_xM/\RR(-X(x))$.

\begin{definition}
\label{defxframing}
An {\em $X$--framing\/} of $L$ is an element of $S(NL,(-X)^{\perp})$.
\end{definition}

Any section $Y$ of $UM$ transverse to $X$ yields an $X$--framing $$\sigma(Y,X) \in S(NL_{Y=-X},(-X)^{\perp})$$ of $L_{Y=-X}$ that is naturally induced by the restriction to $L_{Y=-X}$ of the tangent map to $Y \colon M \rightarrow UM$.
(The tangent map to $Y$ at $x \in L_{Y=-X}$ maps $T_xM$ to $T_xM \oplus U_xM$. Since its composition with the projection onto $U_xM$ maps $T_x L_{Y=-X}$ to $\RR X(x)$, this composition induces a map from $N_x L_{Y=-X}$ to $T_xM/\RR(-X(x))$. This map is orientation-preserving and it is therefore homotopic to a unique $X$-framing.)

\begin{definition}
\label{defcxlsigma}
Conversely, a section $X$ of $UM$ and a link $L$ of $M$ equipped with an $X$--framing $\sigma$ induce a section $C(X,L,\sigma)$ of $UM$ defined as follows (up to homotopy).
Let $N(L)$ be a tubular neighborhood of $L$. 
A fiber $D^2$ of $N(L)$ at $x$ is seen as $\{ uv; u \in [0,1], v \in N_xL\}$.
Let $[-X(x),X(x)]_{\sigma(v)}$ denote the geodesic arc of $U_xM \cong S^2$ from $(-X(x))$ to $X(x)$ through $\sigma(v) \in (-X(x))^{\perp}$.
Then $[-X(x),C(X,L,\sigma)(uv)]$ is the subarc of $[-X(x),X(x)]_{\sigma(v)}$ of length $u \pi$ starting at $-X(x)$. This defines $C(X,L,\sigma)$ on $N(L)$, and $C(X,L,\sigma)$ coincides with $X$ outside $N(L)$.
\end{definition}

\begin{lemma}
\label{lemtubnei}
Let $X$ and $Y$ be two transverse sections of $UM$.
Then $Y$ is homotopic to $C(X,L_{Y=-X},\sigma(Y,X))$.
Furthermore, the Spin$^c$-structure of\/ $Y$ is determined by $[X]^c$ and $L_{Y=-X}$.
\end{lemma}
\bp Outside $L_{Y=-X}$, there is a homotopy from $Y$ to $X$.
When $Y(m) \neq - X(m)$, there is a unique geodesic arc $[Y(m),X(m)]$ with length $(\ell \in [0, \pi[)$ from $Y(m)$ to $X(m)$.
For $t \in [0,1]$, let $Y_t(m) \in [Y(m),X(m)]$ be such that
the length of $[Y(m)=Y_0(m),Y_t(m)]$ is $t\ell$.
Let $D^2$ be the unit disk of $\CC$. Write $N(L_{Y=-X})$ as $D^2 \times L_{Y=-X}$, and let $\chi$ be a smooth increasing bijective function from $[0,1]$ to $[0,1]$ whose derivatives vanish at $0$ and $1$. Set $\tilde{Y}_t(m)=Y_t(m)$ if $m \notin N(L_{Y=-X})$ and $\tilde{Y}_t(v\in D^2,\ell \in L_{Y=-X})=\left\{\begin{array}{ll} Y_{\chi(|v|)t}(v,\ell)\;& \mbox{if }\; v \neq 0\\ 
Y(0,\ell)& \mbox{if }\; v = 0.\end{array}
\right.$\\
Then $\tilde{Y}_1$ is homotopic to $Y$, and when $N(L)$ is small enough, it is easy to see that $\tilde{Y}_1$ is homotopic to $C(X,L_{Y=-X},\sigma(Y,X))$, too.

Let us prove that $[Y]^c=[C(X,L_{Y=-X},\sigma(Y,X))]^c$ does not depend on the $X$-framing $\sigma(Y,X)$ of $L_{Y=-X}$. Two representatives $\sigma_1$ and $\sigma_2$ of any two $X$--framings of a link may be assumed to coincide over the link except over one little interval for each link component. Thus, the associated $C(X,L_{Y=-X},\sigma_1)$ and $C(X,L_{Y=-X},\sigma_2)$ coincide outside a finite union of balls that embeds in a larger ball. Then $[Y]^c$ is determined by $X$ and $L_{Y=-X}$.
Now, changing $X$ inside its homotopy class or changing $X$ over a ball does not affect $[Y]^c$.
\eop

Let $(-X)^{\perp}$ also denote the pull-back of $(-X)^{\perp}$ under the natural projection from $[0,1]\times M$ to $M$.
Let $\Sigma$ be a properly embedded surface in $[0,1] \times M$. Let $S(N\Sigma,(-X)^{\perp})$ denote the space of homotopy classes of sections of the $S^1$--bundle over $\Sigma$ whose fiber over $x$ is the space of orientation-preserving linear isometries from the fiber $N_x\Sigma=T_x([0,1]\times M)/T_x\Sigma$ of $N\Sigma$ to $(-X(x))^{\perp}$.
An {\em $X$--framing\/} of $\Sigma$ is an element of $S(N\Sigma,(-X)^{\perp})$.

Two $X$--framed links $L$ and $L^{\prime}$ are {\em $X$--framed cobordant\/} if and only if there exists an $X$--framed cobordism $\Sigma$ (that is a cobordism equipped with an $X$-framing) properly embedded in $[0,1] \times M$,
from $\{0\} \times L$ to $\{1\} \times L^{\prime}$ that induces the $X$--framings of $L$ and $L^{\prime}$.

\begin{theorem}
\label{thmpontcomb}
Let $X$ be a section of $UM$. Two sections $Y$ and $Z$ of $UM$ transverse to $X$ are homotopic if and only if $(L_{Y=-X},\sigma(Y,X))$ and $(L_{Z=-X},\sigma(Z,X))$ are $X$-framed cobordant.
\end{theorem}
\bp View a homotopy $Y_t$ from $Y=Y_0$ to $Z=Y_1$ as a section $Y_t$ of the pull-back of $UM$ under the natural projection from $[0,1]\times M$ to $M$, and assume without loss that
$(Y_t,-X)$ is transverse to the diagonal of $S^2 \times S^2$ (with respect to some parallelization). Then the preimage $\Sigma$ of the diagonal is a cobordism from $L_{Y=-X}$ and $L_{Z=-X}$ that is canonically $X$--framed by an $X$-framing that induces those of $L_{Y=-X}$ and $L_{Z=-X}$.

Conversely, an $X$--framed cobordism $\Sigma$ from $(L_{Y=-X},\sigma(Y,X))$ to $(L_{Z=-X},\sigma(Z,X))$ induces a section $C(X,\Sigma)$ of the pull-back of $UM$ under the natural projection from $[0,1]\times M$ to $M$ that is defined as $C(X,L,\sigma)$ in Definition~\ref{defcxlsigma} so that the restriction $C_t$ of $C(X,\Sigma)$ on $\{t\} \times M$ defines a homotopy from 
$D_0=C(X,L_{Y=-X},\sigma(Y,X))$ to $D_1=C(X,L_{Z=-X},\sigma(Z,X))$, and, according to Lemma~\ref{lemtubnei}, $Y$ and $Z$ are homotopic.
\eop

\begin{remark}
In \cite[1.4]{laudenbach}, Fran\c{c}ois Laudenbach proves a similar result for nowhere zero sections of a cotangent bundle of a manifold of arbitrary dimension. This result can easily be adapted to any other real bundle over a manifold of the same dimension. Again, I thank Patrick Massot for pointing out this reference to me.
\end{remark}

\begin{corollary}
\label{corpontcomb}
Let $X$ be a section of $UM$. The Spin$^c$--structure of a section $Y$ of $UM$ transverse to $X$ is determined by $[X]^c$ and by the homology class $[L_{Y=-X}]$ of $L_{Y=-X}$ in $H_1(M;\ZZ)$.
\end{corollary}
\eop

A {\em framing\/} of a link $L$ of $M$ is a homotopy class of sections of the unit normal bundle to $L$.
Pushing $L$ in the direction of such a section yields a {\em parallel\/} $L_{\parallel}$ of $L$ up to isotopy of $L_{\parallel}$ in $N(L) \setminus L$, where $N(L)$ is a tubular neighborhood of $L$. This isotopy class of parallels induced by the framing determines the framing. Thus, a {\em framing\/} of $L$ is such an isotopy class of parallels of $L$. 

\begin{remark}
\label{rkxfram}
Let $L$ be a link $X$--framed by some $[\sigma] \in S(NL,(-X)^{\perp})$ represented by $\sigma \colon N(L) \rightarrow (-X)^{\perp}$. (See Definition~\ref{defxframing}.)
Let $\sigma_N$ be a unit section of $N(L)$ that induces a parallel $L_{\parallel}$ of $L$, up to isotopy.
Set $Z(\sigma,\sigma_N)(x)=\sigma(x)(\sigma_N(x))$.
Then $Z(\sigma,\sigma_N)$ is a section of $(-X)^{\perp}$.
Note that $[\sigma]$ is determined by the homotopy classes of $\sigma_N$ and $Z(\sigma,\sigma_N)$, where the homotopy class of $\sigma_N$ may be replaced by the isotopy class of $L_{\parallel}$.
Thus elements of $S(NL,(-X)^{\perp})$ can be thought of as pairs $(L_{\parallel},Z(\sigma,\sigma_N))$ up to simultaneous twists of $L_{\parallel}$ and $Z(\sigma,\sigma_N)$.
\end{remark}

A parallelization $\tau$ with $X$ as first vector identifies $X$-framings of links with framings of links as follows: The second vector $X_2$ of $\tau$ is a section of $(-X)^{\perp}$, and $\tau$ identifies an $X$--framing $[\sigma] \in S(NL,(-X)^{\perp})$ represented by $\sigma$ with the isotopy class of parallels $L_{\parallel}$ of $L$ induced by the section $\sigma^{-1}(X_2)$. Set $$C(\tau,L,L_{\parallel})=C(X,L,[\sigma]).$$
A {\em framed cobordism\/} from $(L,L_{\parallel})$ to $(L^{\prime},L^{\prime}_{\parallel})$
is a cobordism $\Sigma$ from $\{0\} \times L$ to $\{1\} \times L^{\prime}$ in $[0,1] \times M$ equipped with a unit normal section to $T\Sigma$ in $T([0,1] \times M)$, up to homotopy, that induces the given framings of $L$
and $L^{\prime}$. Two framed links are {\em framed cobordant\/} if and only if their exists a framed cobordism from one to the other one.

As above, a parallelization $\tau$ with $X$ as first vector identifies $X$-framings of cobordisms to framings of cobordisms.

This allows us to state the following Pontrjagin theorem \cite[Section 7, Theorem B]{Mil} as a corollary of Theorem~\ref{thmpontcomb}.

\begin{theorem}[Pontrjagin construction]
\label{thmpont}
Let $\tau$ be a parallelization of $M$.
Any section of $UM$ is homotopic to $C(\tau,L,L_{\parallel})$ for a framed link $(L,L_{\parallel})$ of the interior of $M$.
Two sections $C(\tau,L,L_{\parallel})$ and $C(\tau,L^{\prime},L^{\prime}_{\parallel})$ are homotopic if and only if $(L,L_{\parallel})$ and $(L^{\prime},L^{\prime}_{\parallel})$ are framed cobordant.
\end{theorem}
\eop

Pontrjagin proved generalizations of this theorem to every dimension. See \cite[Section 7]{Mil}.

Let $\Sigma_M$ be an embedded cobordism from a link $L$ to a link $L_1$ in $M$.
The graph of a Morse function $f$ from $\Sigma_M$ to $[0,1]$ such that $f^{-1}(0)=L$ and $f^{-1}(1)=L_1$ yields a proper embedding $\Sigma$ of $\Sigma_M$ into $[0,1]\times M$. The positive normal to $\Sigma_M$ in $M$ at $m$ seen in $T_{(f(m),m)}\{f(m)\}\times M$ frames $\Sigma$.
This framing of $\Sigma$ identifies the $X$-framings $\Sigma$ with homotopy classes of sections of $(-X)^{\perp}$ over $\Sigma$.
When $\Sigma_M$ is connected, and when $K$ is a boundary component of $\Sigma$, any $X$-framing defined on $\partial \Sigma \setminus K$ extends as an $X$-framing of $\Sigma$, and the extension of the $X$-framing over $K$ is determined by the restriction of the $X$-framing to $\partial \Sigma \setminus K$.

Embed a sphere $S$ with three holes in $M$, the $3$ boundary components of $S$ are $3$ knots 
$K_1$, $K_2$ and $-K_1 \sharp_b K_2$ of $M$ that are framed by the embedding of $S$. 

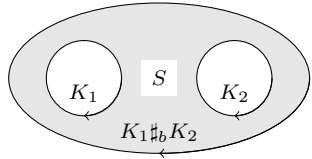
\begin{center}
\begin{tikzpicture}
\useasboundingbox (0,0) rectangle (4,2);
\draw [fill=gray!20] (2,1) ellipse (2 and 1);
\draw (2,1) node[rectangle,fill=white]{\scriptsize $S$};
\draw [->] (4,1) arc (0:-90:2 and 1);
\draw (2,0) node[above]{\scriptsize $K_1 \sharp_b K_2$};
\draw [fill=white] (1,1) circle (.5) (3,1) circle (.5);
\draw [->] (1.5,1) arc (0:-90:.5);
\draw [->] (3.5,1) arc (0:-90:.5);
\draw (1,.8) node{\scriptsize $K_1$} (3,.8) node{\scriptsize $K_2$};
\end{tikzpicture}
\end{center}

Then $K_1 \sharp_b K_2$ is a {\em framed band sum} of $K_1$ and $K_2$, it is framed cobordant to the union of $K_1$ and $K_2$.
Note that any $X$--framed link is $X$--framed cobordant to an $X$--framed knot by such band sums.
Similarly, any framed link is framed cobordant to a framed knot.

\begin{lemma}
 \label{lemframedZ}
Two framed links $(L,L_{\parallel})$ and  $(L^{\prime},L_{\parallel}^{\prime})$ in a $\ZZ$-sphere or in a $\ZZ$-ball are framed cobordant
if and only if $lk(L,L_{\parallel})=lk(L^{\prime},L_{\parallel}^{\prime})$.
\end{lemma}
\bp When the framed links are framed cobordant, $lk(L,L_{\parallel})=lk(L^{\prime},L_{\parallel}^{\prime})$, since $lk(L,L_{\parallel})$ is the algebraic intersection of two $2$-chains bounded by $L \times \{0\}$ and $L_{\parallel} \times \{0\}$ in $[-1,0] \times M$. Conversely, let $(L,L_{\parallel})$ and  $(L^{\prime},L_{\parallel}^{\prime})$ be two framed links such that $lk(L,L_{\parallel})=lk(L^{\prime},L_{\parallel}^{\prime})$. They are respectively framed cobordant to framed knots $(K,K_{\parallel})$ and $(K^{\prime},K_{\parallel}^{\prime})$ such that $lk(K,K_{\parallel})=lk(L,L_{\parallel})$ and $lk(K^{\prime},K_{\parallel}^{\prime})=lk(L^{\prime},L_{\parallel}^{\prime})$, so that $lk(K,K_{\parallel})=lk(K^{\prime},K_{\parallel}^{\prime})$. There is a connected cobordism from $K$ to $K^{\prime}$ that may be equipped with a framing that extends the framing induced by $K_{\parallel}$, and that therefore induces a framing of $K^{\prime}$ corresponding to a parallel $K^{\prime}_1$ of $K^{\prime}$ such that $lk(K,K_{\parallel})=lk(K^{\prime},K^{\prime}_1)$. Thus $lk(K^{\prime},K^{\prime}_1)=lk(K^{\prime},K^{\prime}_{\parallel})$ and $K^{\prime}_1$ is isotopic to $K^{\prime}_{\parallel}$, so that $(K^{\prime},K_{\parallel}^{\prime})$ is framed cobordant to  $(K,K_{\parallel})$.
\eop

\subsection{More details about the introductions}
\label{subdetint}

Let us finish justifying the claims of the introductions.

\begin{lemma}
\label{lemadL}
For any two transverse sections $X$ and\/ $Y$ of\/ $UM$, $L_{Y=-X}=-L_{X=-Y}$.
For three pairwise transverse sections $X$, $Y$ and $Z$ of $UM$, $[L_{Z=-X}]=[L_{Z=-Y}]+[L_{Y=-X}]$ in $H_1(M;\ZZ)$.
\end{lemma}
\bp
For two sections $X$ and $Z$ of $UM$, transverse to $Y$, up to homotopy, we can assume that $L_{X=-Y}$ and $L_{Z=-Y}$ are disjoint, and pick disjoint tubular neighborhoods $N(L_{X=-Y})$  and $N(L_{Z=-Y})$ of $L_{X=-Y}$ and $L_{Z=-Y}$, respectively.
Then, according to Lemmas~\ref{lemLXYwd} and \ref{lemtubnei} we can assume that $Z=C(Y,L_{Z=-Y},\sigma(Z,Y))$ and that $X=C(Y,L_{X=-Y},\sigma(X,Y))$ so that $Z=Y$ outside $N(L_{Z=-Y})$ and $X=Y$ outside $N(L_{X=-Y})$.
Then $L_{Z=-X}=L_{Z=-Y} \coprod L_{Y=-X}$.
\eop

\begin{lemma}
 There is a canonical free transitive action of $H_1(M;\ZZ)$ on the set $\CS(M)$ of Spin$^c$-structures of $M$ such that for any two transverse sections $Y$ and $Z$ of $UM$, $$[L_{Z=-Y}][Y]^c=[Z]^c.$$
\end{lemma}
\bp Let $Y$ be a section of $UM$ and let $[K] \in H_1(M;\ZZ)$.
Represent $[K]$ by a knot $K$ and equip $K$ with an arbitrary $Y$-framing $\sigma$.
Define $[K][Y]^c$ as $[Z]^c$ with $Z=C(Y,K,\sigma)$. According to Definition~\ref{defcxlsigma}, $K=L_{Z=-Y}$ and, according to Corollary~\ref{corpontcomb}, $[Z]^c$ is determined by $[Y]^c$ and $[K]$.
According to Lemma~\ref{lemLXYwd}, if $[K][Y]^c=[Z]^c$, then $K$ is homologous to $L_{Z=-Y}$.
Lemma~\ref{lemadL} ensures that this defines an action of $H_1(M;\ZZ)$. This action is obviously transitive since
$[Z]^c=[L_{Z=-Y}][Y]^c$ and it is free.
\eop

\begin{corollary}
This action equips $\CS(M)$ with an affine structure with translation group $H_1(M;\ZZ)$.
With respect to this structure, for any two transverse sections $X$ and $Y$ of $UM$, $$[Y]^c-[X]^c=[L_{Y=-X}].$$
\end{corollary}

Classically, $\CS(M)$ is rather equipped with an affine structure with translation group $H^2(M,\partial M;\ZZ)$, and $([Y]^c-[X]^c)_2 \in H^2(M,\partial M;\pi_2(S^2)=\ZZ)$ is the obstruction to homotoping a section $Y$ of $UM$ to another such $X$ over a two-skeleton of $M$.

Below, we confirm that the two structures are naturally related by the Poincar\'e duality isomorphism
$P \colon H^2(M,\partial M;\ZZ) \rightarrow H_1(M;\ZZ)$.

\begin{lemma}
For two transverse sections $X$ and $Y$ of $UM$, $$P(([Y]^c-[X]^c)_2)=[Y]^c-[X]^c=[L_{Y=-X}].$$
\end{lemma}
\bp
Up to homotopy, assume $Y=C(X,L_{Y=-X},\sigma(Y,X))$ as in Lemma~\ref{lemtubnei}.
Let $S$ be a $2$--chain transverse to $L_{Y=-X}$. We may assume that $X$ and $Y$ coincide outside open disks around $S \cap L_{Y=-X}$. Extend $X$ to a parallelization on the closure of these disks, and see $Y$ as a map from $D^2/\partial D^2$ to $S^2$ on each of these disks. The sum of the degrees of these maps is the algebraic intersection of $L_{Y=-X}$ and $S$.
By definition, this is also the evaluation of a cochain that represents $([Y]^c-[X]^c)_2 \in H^2(M, \partial M;\ZZ)$ on $S$. This shows that $L_{Y=-X}$ is Poincar\'e dual to $([Y]^c-[X]^c)_2$.
\eop

The {\em Euler class\/} $e(X^{\perp})$ is the obstruction to the existence of a nowhere zero section of $X^{\perp}$. It lives in $H^2(M;\ZZ)$.
In particular, $X$ extends as a parallelization if and only if $e(X^{\perp})=0$.
We shall not give a more precise definition for the standard Euler class, since Lemmas~\ref{lemdefeul} below can be used as a definition in our cases.

\begin{lemma}
\label{lemdefeul}
 Let $X$ and $Y$ be two homotopic transverse sections of $UM$, then $L_{Y=X}$ is Poincar\'e dual to $e(X^{\perp})$. Therefore, $P(e(X^{\perp}))=[X]^c-[-X]^c$.
\end{lemma}
\bp
For a section of $X^{\perp}$, $X$ may be pushed slightly in the direction of the section. If $Y$ denotes the obtained combing, then $L_{Y=X}$ is the vanishing locus of the section that is Poincar\'e dual to $e(X^{\perp})$.
\eop

\begin{lemma}
\label{lemeultor}
Let $X$ and $Y$ be two transverse sections of $UM$,
$$2[L_{X=Y}]=P(e(X^{\perp})+e(Y^{\perp}))$$
and $2[L_{X=-Y}]=P(e(X^{\perp})-e(Y^{\perp}))$.
In particular, for two transverse torsion sections $X$ and $Y$ of $UM$, $L_{X=Y}$ and $L_{X=-Y}$ represent torsion elements in $H_1(M;\ZZ)$.
\end{lemma}
\bp
$[L_{X=Y}]=[X]^c-[-Y]^c=[Y]^c-[-X]^c$ so that $2([L_{X=Y}])=[X]^c-[-X]^c+[Y]^c-[-Y]^c=P(e(X^{\perp})+e(Y^{\perp})).$
\eop
\begin{lemma}
\label{lemlkind}
Let $X$ and $Y$ be two transverse torsion sections of $UM$, then $lk(L_{X=Y},L_{X=-Y})$ only depends on the homotopy classes of $X$ and $Y$.
\end{lemma}
\bp Fix a trivialization of $UM$ so that sections become functions from $M$ to $S^2$. Let us prove that $lk(L_{X=Y},L_{X=-Y})$ does not vary under a generic homotopy of $X$.
Such a homotopy induces two homotopies $h_+$ and $h_-$ from $[0,1]\times M$ to $S^2 \times S^2$ where $h_{\pm}(t,m)=(X_t(m),\pm Y(m))$. Without loss, assume that $h_+$ and $h_-$ are transverse to the diagonal.
There exists a finite sequence $0=t_0<t_1<t_2<\dots<t_k=1$ of times such that the projections on $M$ of the preimages of the diagonal
under $h_{+|[t_i,t_{i+1}] \times M}$ and $h_{-|[t_i,t_{i+1}] \times M}$ are disjoint so that they yield two disjoint cobordisms in $M$, one from $L_{X_{t_i}=Y}$ to $L_{X_{t_{i+1}}=Y}$, and the other one from $L_{X_{t_i}=-Y}$ to $L_{X_{t_{i+1}}=-Y}$ showing that $lk(L_{X_{t_i}=Y},L_{X_{t_i}=-Y})=lk(L_{X_{t_{i+1}}=Y},L_{X_{t_{i+1}}=-Y})$.
\eop

\begin{lemma}
\label{lemdetcom}
Let $X$ be a section of $UM$ that extends as a parallelization $\tau$. The homotopy class of a torsion section $Y$ transverse to $X$ is determined by $X$, by the homology class $[L_{Y=-X}]$ of $L_{Y=-X}$ in $H_1(M;\ZZ)$, and by the linking number $lk(L_{Y=-X},L_{Y=X})$.
\end{lemma}
\bp After a homotopy, $Y$ reads $C(\tau,L_{Y=-X},L_{Y=X_2})$ where $X_2$ is the second vector of $\tau$, and, $L_{Y=X}$ and $L_{Y=X_2}$ are parallel knots as in Theorem~\ref{thmpont}. According to Theorem~\ref{thmpont}, the combing $[Y]$ is determined by the framed cobordism class of $L_{Y=-X}$, that is determined by $[L_{Y=-X}]$ and by $lk(L_{Y=-X},L_{Y=X_2})$ since $L_{Y=-X}$ is rationally null-homologous. After another homotopy that makes $Y$ transverse to $X_2$ and $X$, $lk(L_{Y=-X},L_{Y=X_2})=lk(L_{Y=-X},L_{Y=X})$.
\eop

\subsection{Action of $\pi_3(S^2)$ on combings}
\label{subpitrois}

\begin{notation}
\label{notrho}
See $B^3$ as the quotient of $[0,2\pi] \times S^2$ where the quotient map identifies $\{0\} \times S^2$ with a point. Then the map from $B^3$ to the group $SO(3)$ of orientation-preserving linear isometries of $\RR^3$ that maps $(\theta \in [0,2\pi],x \in S^2)$ to the rotation $\rho(\theta,x)$ with axis directed by $x$ and with angle $\theta$ is denoted by $\rho$.
It induces the standard double covering map $\tilde{\rho}$ from $S^3=B^3/\partial B^3$ to $SO(3)$ that orients $SO(3)$.
\end{notation}

The image of the first basis vector $p_{S^2} \colon SO(3) \rightarrow S^2$ induces an isomorphism from $\pi_3(SO(3))=\ZZ[\tilde{\rho}]$ to $\pi_3(S^2)$. Let $\gamma$ be the image of $[\tilde{\rho}]$ under this isomorphism.
Let $X$ be a combing. Extend $X$ to a parallelization $(X,Y,Z)$ on a $3$-ball $B$ identified with $B^3$, and see $\rho$ as a map $\rho \colon (B,\partial B) \rightarrow (SO(X,Y,Z),\mbox{Id})$. Define $\gamma^k X$ as the section that coincides with $X$ outside $B$ and such that, for any $m \in B$,
$$\gamma^k X(m)=(\rho(m))^k(X)$$ on $B$. Note that $[\gamma^k X]$ is independent of the chosen parallelization. Since $M$ is connected, any two small enough balls may be put inside a bigger one and $[\gamma^k X]$ is independent of $B$. Set $\gamma^k[X]=[\gamma^k X]$. Note that $\gamma^{k+k^{\prime}}[ X]=\gamma^{k} (\gamma^{k^{\prime}}[X])$.
Let $X$ and $Y$ be two sections of $UM$ that are homotopic except over a $3$-ball $B^3$. Up to homotopy, we may assume that they are identical outside $B^3$. On $B^3$, $X$ extends to a parallelization and $Y$ reads as a map from $(B^3,\partial B^3)$ to $(S^2,X)$. It therefore defines an element $\gamma^k$ of $\pi_3(S^2)$, and $[Y]=\gamma^k[ X].$
Thus, $\pi_3(S^2)$ acts transitively on the combings that represent a given Spin$^c$-structure. In particular it acts transitively on the combings of a $\ZZ$-sphere.

A \emph{positive (or oriented) meridian} of some knot $K$ in $M$ is the boundary of a disk that intersects $K$ once with positive sign.

\begin{lemma}
\label{lemhopfneg}
Let $\tau$ be a parallelization of $M$ and let $[X(\tau)]$ denote the induced combing. Let $(U,U_-)$ be the negative Hopf link \neghopf ($lk(U,U_-)=-1$).
Then, with the notation before Theorem~\ref{thmpont}, $[\gamma X(\tau)]=[C(\tau,U,U_-)]$.
\end{lemma}
\bp
First note that $[C(\tau,U,U_-)]$ reads $[\gamma^k X(\tau)]$ for an integer $k$ that does not depend on $(M,\tau)$. We prove $k=1$ when $M=B^3$, when $\tau$ is the standard parallelization, and when $X=X(\tau)$ is the constant upward vector field, with the help of Lemma~\ref{lemdetcom}, by showing that
$$lk(L_{\gamma X(\tau)=X^{\prime}},L_{\gamma X(\tau)=-X^{\prime}})=lk(U,U_-)=-1.$$ 
for some constant field $X^{\prime}$ near $X$.
Let $N$ be the North pole of $S^2$, $(p_{S^2}\circ \rho)^{-1}(N)$ intersects the interior of $B^3$ as the vertical axis oriented from South to North while $(p_{S^2}\circ \rho)^{-1}(-N)$ intersects $B^3$ as $\pi \times (-E)$, where $E$ is the equator oriented as a positive meridian
 of $(p_{S^2}\circ \rho)^{-1}(N)$.
Then for $N^{\prime}$ near $N$, $lk((p_{S^2}\circ \rho)^{-1}(N^{\prime}),(p_{S^2}\circ \rho)^{-1}(-N^{\prime}))=-1$.
\eop

\begin{corollary}
\label{coractpito}
Let $\tau$ be a parallelization of $M$, let $(L,L_{\parallel})$ be a framed link of $L$, let $(U,U_-)$ be the negative Hopf link in a ball of $M$ disjoint from $L$, and let $(U,U^+)$ be the positive Hopf link \poshopf in a ball of $M$ disjoint from $L$.
Then $[\gamma C(\tau,L,L_{\parallel})]=[C(\tau,L \cup U,L_{\parallel} \cup U_-)]$ and $[\gamma^{-1} C(\tau,L,L_{\parallel})]=[C(\tau,L \cup U,L_{\parallel} \cup U_+)]$.

If $L$ is non-empty, let $L_{\parallel,-1}$ (resp. $L_{\parallel,+1}$) be a parallel of $L$ obtained from $L_{\parallel}$ by adding a negative (resp. positive) meridian of $L$, homologically in $N(L)\setminus L$, then $[C(\tau,L \cup U,L_{\parallel} \cup U_-)]=[C(\tau,L,L_{\parallel,-1})]$ and $[C(\tau,L \cup U,L_{\parallel} \cup U_+)]=[C(\tau,L,L_{\parallel,+1})]$.
\end{corollary}
\bp
Note that $(L,L_{\parallel,\pm 1})$ is framed cobordant to $(L \cup U,L_{\parallel} \cup U_{\pm})$ by band sum. Thus, 
the second formula can be deduced from the fact that the disjoint union of two oppositely framed unknots is framed cobordant to the empty link.
\eop

\begin{corollary}
\label{coractpit}
Let $X$ be a torsion section of $UM$, let $k \in \ZZ$ and let $Y$ be a section of $UM$ that represents $[\gamma^{k} X]$.
Then $lk(L_{Y=X},L_{Y=-X})=-k$.
\end{corollary}
\bp We already know that the linking number $lk(L_{Y=X},L_{Y=-X})$ does not depend on the transverse representatives of $[X]$ and $[Y]$. Furthermore, by Theorem~\ref{thmpont}, $[X]$ can be represented as $C(\tau,L,L_{\parallel})$ as in Corollary~\ref{coractpito}. Assume $k \neq 0$. Let $(\cup_{i=1}^{|k|} U^{(i)}, \cup_{i=1}^{|k|} U^{(i)}_{\varepsilon})$ denote the union of $|k|$ Hopf links with sign $\varepsilon=-k/|k|$ contained in disjoint balls $B_i$, for $i=1, \dots, k$.
Let $Y$ be obtained from $C(\tau,L \cup \cup_{i=1}^{|k|} U^{(i)},L_{\parallel} \cup \cup_{i=1}^{|k|} U^{(i)}_{\varepsilon})$ by a small perturbation, induced by the parallelization $\tau$ outside $N(L  \cup \cup_{i=1}^{|k|} U^{(i)})$ so that it is transverse to $X$, very close to $X$, and distinct from $\pm X$ outside $N(L \cup (\cup_{i=1}^{|k|} U^{(i)}))$. Then $L_{Y=-X}$ is a parallel of $\cup_{i=1}^{|k|} U^{(i)}$ and 
$$lk(L_{Y=X},L_{Y=-X})=\sum_{i=1}^{|k|}lk(L_{Y=X} \cap B_i,L_{Y=-X} \cap B_i)=\sum_{i=1}^{|k|}lk(U^{(i)},U^{(i)}_{\varepsilon})=-k.$$
\eop

\begin{proposition}
\label{propSpinstruc}
 Let $[X]^c$ be a Spin$^c$ structure. Then the set of combings that belong to $[X]^c$ is an affine space over $\ZZ/\langle e(X^{\perp}),H_2(M;\ZZ)\rangle$, where the translation by the class of $1$ is the action of $\gamma$.
\end{proposition}
\bp
Again, fix a parallelization $\tau$ of $M$, and an induced combing $Y$. This identifies the set
$\CS(M)$ of Spin$^c$ structures with $H_1(M;\ZZ)$ by mapping $[X]^c$ to the homology class $[L_{X=Y}]$.
Any framed link is framed cobordant to a framed knot. According to the Pontrjagin characterization of the combings (Theorem~\ref{thmpont}), the combings that belong to the Spin$^c$ structure $\xi(\tau,[K])$ corresponding to a given class $[K]$ of $H_1(M;\ZZ)$ is the set of framed links homologous to $[K]$ up to framed cobordism. Let $K$ be a knot that represents $[K]$, then all framed links homologous to $[K]$ are framed cobordant to $K$ equipped with some framing, and the combings of $\xi(\tau,[K])$ are the equivalence classes of framings of $K$ up to framed cobordism.

For two parallels $K^{\prime}$ and $K^{\prime\prime}$ of $K$ on the boundary $\partial N(K)$ of a tubular neighborhood $N(K)$ of $K$, the homology class of $K^{\prime\prime} - K^{\prime}$ in $\partial N(K)$ reads $lk_{N(K)}(K^{\prime\prime} - K^{\prime},K)m(K)$ where $m(K)$ is the oriented meridian of $K$. The integer $lk_{N(K)}(K^{\prime\prime} - K^{\prime},K)$ measures the difference of the framings induced by $K^{\prime}$ and $K^{\prime\prime}$.

When $[K]$ is a torsion element of $H_1(M;\ZZ)$, the self-linking number $lk(K^{\prime},K)$ makes sense, and it is a complete invariant of framings of $K$, up to framed cobordism. This shows that the action of $\pi_3(S^2)$ on the set of combings in a torsion Spin$^c$-structure is free, and that this set is an affine space over $\ZZ$.

In general, let $B$ be a cobordism from $0 \times K^{\prime}$ to $1 \times K^{\prime\prime}$ in $[0,1] \times N(K)$.
Then $lk_{N(K)}(K^{\prime\prime} - K^{\prime},K)=\langle [0,1] \times K,B\rangle_{[0,1] \times M}$.
Let $C$ be a framed cobordism from $0 \times K$ to $1 \times K$ in $[0,1] \times M$, and let $C^{\prime}$ be obtained from $C$ by pushing $C$ in the direction of the framing.
Assume that $\partial C^{\prime}=1 \times K^{\prime\prime} - 0 \times K^{\prime}$ so that
$C$ is a framed cobordism from $(K,K^{\prime})$ to $(K,K^{\prime\prime})$ and
$$0=\langle C,C^{\prime}\rangle_{[0,1] \times M}=\langle [0,1] \times K +(C-[0,1] \times K),B+(C^{\prime}-B)\rangle_{[0,1] \times M}.$$
Since $(C-[0,1] \times K)$ and $(C^{\prime}-B)$ are $2$-cycles in $[0,1] \times M$, $\langle (C-[0,1] \times K),(C^{\prime}-B)\rangle_{[0,1] \times M}=0$, and since they are homologous
$\langle [0,1] \times K,(C^{\prime}-B)\rangle_{[0,1] \times M}=\langle (C-[0,1] \times K),B\rangle_{[0,1] \times M}$, so that
$$lk_{N(K)}(K^{\prime\prime} - K^{\prime},K)=-2\langle [0,1] \times K ,(C^{\prime}-B)\rangle_{[0,1] \times M}.$$
In particular, the framing difference induced by $C$ only depends on the homology class of the projection $S$ of $C$ in $M$, and it is $-2\langle K,S\rangle_M$.
Thus if the framings induced by $K^{\prime}$ and $K^{\prime\prime}$ are framed cobordant, $lk_{N(K)}(K^{\prime\prime} - K^{\prime},K)$ is in $\langle 2K ,H_2(M;\ZZ)\rangle_M$. Conversely, for any class $S$ of $H_2(M;\ZZ)$, there exists an embedded connected cobordism $C$ that projects on $S$. Any framing on $0\times K$ can be extended to $C$, and it induces a framing on $1\times K$, such that the framing difference is $-2\langle K,S\rangle_M$.
Since the Euler class of $\xi(\tau,[K])$ is Poincar\'e dual to $2[K]$, the conclusion follows.
\eop

\section{Towards the variation formula}
\setcounter{equation}{0}
\label{seckeyprop}

\subsection{The key proposition}
In this subsection that will be useful in our study of the invariant $\Theta$ in Section~\ref{sectheta}, we prove the following proposition that is the key to the extension of the map $p_1$ in Section~\ref{secpone}.

\begin{proposition}
\label{proptrans}
 Let $X$, $Y$ and $Z$ be three pairwise transverse torsion sections of $UM$, 
$$lk(L_{X=Y},L_{X=-Y})+lk(L_{Y=Z},L_{Y=-Z})=lk(L_{X=Z},L_{X=-Z}).$$
\end{proposition}

Consider the $6$-manifold $[0,1]\times UM$. Recall that $UM$ is homeomorphic to $M \times S^2$.
Let $(S_i)_{i=1,\dots,\beta_1(M)}$ be $\beta_1(M)$ surfaces in the interior of $M$ that represent a basis of $H_2(M;\QQ)$.
For a section $Z$ of $UM$, let $Z(S_i)$ denote the image in $UM$ of the graph of the restriction of $Z$ to $S_i$. Let $[S]$ denote the homology class of the fiber of $UM$ in $H_2(UM;\QQ)$, oriented as the boundary of a unit ball of $T_xM$.
$$H_2(UM;\QQ)=\QQ[S] \oplus \bigoplus_{i=1}^{\beta_1(M)}\QQ[Z(S_i)].$$
\begin{lemma}
\label{lemsecYZ}
If\/ $Y$ and $Z$ are two transverse sections of $UM$,
then $$[Z(S_i)]-[Y(S_i)]=\langle L_{Z=-Y} ,S_i\rangle_M[S]$$ in $H_2(UM;\QQ)$ (and in $H_2([0,1]\times UM;\QQ)$).
\end{lemma}
\bp
Fix a trivialization of $UM$ so that both $Y$ and $Z$ become functions from $M$ to $S^2$, then $[Z(S_i)]-[Y(S_i)]=(\mbox{deg}(Z_{|S_i})-\mbox{deg}(Y_{|S_i}))[S]$. If $X$ is a section of $UM$ induced by the trivialization, then $\mbox{deg}(Z_{|S_i})=\langle L_{Z=-X},S_i\rangle_M$. Conclude with Lemma~\ref{lemadL}.
\eop

In particular, according to Lemma~\ref{lemeultor}, the subspace $H_T$ of $H_2([0,1]\times UM;\QQ)$ generated by the $[Z(S_i)]$ for torsion combings $Z$ is canonical.
Set $H(M)=H_2([0,1]\times UM;\QQ)/H_T$. Then $H(M)=\QQ[S]$. 

Let $X$ and $Y$ be two sections of $UM$. Let $X(M)$ abusively denote the graph of $X$ in $UM$. 
Let $\partial(X,Y)$ be the following codimension $2$ submanifold of $\partial([0,1]\times UM)$. If $\partial M =\emptyset$, $\partial(X,Y)=\{1\} \times Y(M) - \{0\} \times X(M)$.
If $\partial M =S^2$, let $V(X)$ and $V(Y)$ be the elements of $S^2$ such that $X=V(X)$ and $Y=V(Y)$ on $\partial M$.
Recall that $\tau_s$ identifies $UM_{|\partial M}$ with $S^2 \times \partial M$. Let $P=P(X,Y)$ be a $1$--chain in $[0,1]\times S^2$ such that $\partial P= \{1\} \times V(Y)-\{0\}\times V(X)$.
Then  $\partial(X,Y)=\partial(X,Y,P)=\{1\} \times Y(M) - \{0\} \times X(M)- P \times \partial M$.

\begin{lemma} \label{lemFXY}
For two transverse sections $X$ and $Y$ of $UM$ such that $([Y]^c-[X]^c)$ vanishes in $H_1(M;\QQ)$, $\partial(X,Y)$ is rationally null-homologous in $[0,1]\times UM$.
It bounds a rational chain $F(X,Y)$ transverse to $\partial ([0,1]\times UM)$ that is well-determined, up to the addition of a chain $\Sigma \times \partial M$ for a $2$--chain $\Sigma$ of $[0,1]\times S^2$, up to the addition of a combination of $\{t_i\} \times UM_{|S_i}$ for distinct $t_i$, and up to cobordism.
\end{lemma}
\bp $H_3([0,1]\times UM;\QQ)\cong H_1(M;\QQ) \otimes H_2(S^2;\QQ)$ when $\partial M=S^2$. The direct factor $\QQ[X(M)]$ should be added when $\partial M=\emptyset$.
The class of a $3$--submanifold of $[0,1]\times UM$ vanishes in $H_3([0,1]\times UM;\QQ)$ if and only if its algebraic intersection with the $[0,1]\times Z(S_i)$ vanishes, for all $i$, when $\partial M=S^2$, for some combing $Z$. 
For $\partial(X,Y)$, this algebraic intersection reads
$$\begin{array}{ll}\langle [0,1]\times Z(S_i), \partial(X,Y)\rangle_M &=\langle S_i, L_{Z=Y}-L_{Z=X}\rangle=\langle S_i, [Z]^c-[-Y]^c-([Z]^c-[-X]^c)\rangle\\&=\langle S_i, [Y]^c-[X]^c\rangle=0.\end{array}$$
When $\partial M=\emptyset$, the algebraic intersection with $[0,1]\times UM_{|\{x\}}$ must vanish, too. This is easily verified.
Thus, $\partial(X,Y)$ bounds a rational chain $F(X,Y)$, and since $H_4( UM;\QQ)=\bigoplus_{i=1}^{\beta_1(M)}\QQ[UM_{|S_i}]$, the second assertion follows.
\eop

\begin{lemma}
\label{lemevHM}
 For any two transverse torsion sections $X$ and $Y$ of $UM$, for any two-cycle $C$ of $[0,1] \times UM$,
the class of $C$ in $H(M)$ is $\langle C, F(X,Y)\rangle_{[0,1] \times UM}[S]$ for a $F(X,Y)$ as in Lemma~\ref{lemFXY}.
\end{lemma}
\bp First note that $\langle C, F(X,Y)\rangle_{[0,1] \times UM}[S]$ only depends on the homology class of $C$, for a given $F(X,Y)$, and that 
$\langle [S], F(X,Y)\rangle=1$. Now, it suffices to prove that $\langle [Z(S_i)], F(X,Y)\rangle=0$ for any torsion combing $Z$, and for any $i$. Since $\langle [Z(S_i)], F(X,Y)\rangle=\langle [Z(S_i)], X(M)\rangle_{UM}=\langle [Z(S_i)], Y(M)\rangle_{UM}$, $\langle [Z(S_i)], F(X,Y)\rangle$ does not depend on the torsion combings $X$ and $Y$. In particular, $\langle [Z(S_i)], F(X,Y)\rangle=\langle [Z(S_i)], F(-Z,-Z)\rangle=0$.
\eop

\begin{definition}
\label{defblowup}
In this article, {\em blowing up\/} a submanifold $A$ means replacing it by its unit normal bundle. 
Let $c$ be the codimension of $A$. The total space of the normal bundle to $A$ locally reads $\RR^c \times U$ for an open subspace $U$ of $A$. It embeds into the ambient manifold as a tubular neighborhood of $A$. Its fiber $\RR^c$ reads $\{0\} \cup (]0,\infty[ \times S^{c-1})$ where the unit sphere $S^{c-1}$ of $\RR^c$ is the fiber of the unit normal bundle to $A$. Then the blow-up replaces $(0 \in \RR^c)$ by $S^{c-1}$ so that the blown-up manifold locally reads $([0,\infty[ \times S^{c-1} \times U)$. (In particular, unlike the blow-ups in algebraic geometry, our differential blow-ups create boundaries.)
Topologically, this blow-up amounts to removing an open tubular neighborhood of $A$ (thought of as infinitely small), but the process is canonical, so that the created boundary is the unit normal bundle to $A$ and there is a canonical projection from the blown-up manifold to the initial manifold.

\end{definition}

\begin{proposition}
\label{propFXYcap}
Let $X$ and $Y$ be two transverse torsion sections of $UM$.
For any $F(X,Y)$ and $F(-X,-Y)$ as in Lemma~\ref{lemFXY}, such that the $1$--chains $P(X,Y)$ and $P(-X,-Y)$ are disjoint, the class of $F(X,Y)\cap F(-X,-Y)$ in $H(M)$ is
$$lk(L_{X=Y},L_{X=-Y})[S].$$
\end{proposition}
\bp Let us first prove that
the class of $F(X,Y)\cap F(-X,-Y)$ is well-determined in $H(M)$.
When $F(X,Y)$ is changed to $F(X,Y) + (\Sigma \times \partial M)$ for a two-chain $\Sigma$ of $[0,1]\times S^2$ transverse to $P(-X,-Y)$,  $(\Sigma \times \partial M) \cap F(-X,-Y)$ is added to $F(X,Y)\cap F(-X,-Y)$.
Now, $(\Sigma \times \partial M) \cap F(-X,-Y)$ is a union of $(t_j,V_j) \times \partial M $ that bounds
since the parallelization $\tau_s$ extends to $M$. Thus, the class of $F(X,Y)\cap F(-X,-Y)$ in $H(M)$ in unchanged.
Since the class of $\{t_i\} \times UM_{|S_i} \cap F(-X,-Y)$ is in $H_T$, the class of $F(X,Y)\cap F(-X,-Y)$ is well-determined in $H(M)$.

Now, we construct an explicit $F(X,Y)$ by using the homotopy of Lemma~\ref{lemtubnei} that we recall.
Assume $M$ is Riemannian.
When $X(m) \neq - Y(m)$, there is a unique geodesic arc $[X(m),Y(m)]$ with length $(\ell \in [0, \pi[)$ from $X(m)$ to $Y(m)$.
For $t \in [0,1]$, let $X_t(m) \in [X(m),Y(m)]$ be such that
the length of $[X_0(m)=X(m),X_t(m)]$ is $t\ell$.
This defines $X_t$ on $(M\setminus L_{X=-Y})$.

Observe that this definition naturally extends to the boundary of the manifold $M(L_{X=-Y})$ obtained from $M$ by blowing up $L_{X=-Y}$:
Indeed, $X$ induces an orientation-preserving map from the normal bundle $N_xL_{X=-Y}$ to $L_{X=-Y}$ in $M$ at $x$ to $(-Y(x))^{\perp}$.
Then for a unit element $n$ of $N_xL_{X=-Y}$, $X_t(n)$ describes the half great circle from $X(x)$ to $Y(x)$ through the image of $n$ under the above map.
In particular, the whole sphere is covered with degree $1$ by the image of $([0,1] \times (N_xL_{X=-Y}/\RR^{\ast +}))$.
Let $G_h$ be the closure of 
$\left(\cup_{t \in [0,1]} s_{X_t}\left(M \setminus L_{X=-Y}\right)\right)$.
$$G_h=\cup_{t \in [0,1]}X_{t}(M(L_{X=-Y})).$$
Define the $3$--cycle of $UM$ $$p(\partial(X,Y))=Y(M)-X(M) - [V(X),V(Y)] \times \partial M$$ where
$[V(X),V(Y)]$ is the shortest geodesic path from $V(X)$ to $V(Y)$ in the fiber of $UM$ over $\partial M$ that is identified with $S^2$ by $\tau_s$.
Then
$$\partial G_h - p(\partial(X,Y))=\cup_{t \in [0,1]}X_{t}(-\partial M(L_{X=-Y}))=UM_{|L_{X=-Y}}$$
because it is oriented like $\cup_{t \in [0,1]}X_{t}( \partial N(L_{X=-Y}))$.
Let $\Sigma_{X=-Y}$ be a two-chain transverse to $L_{X=Y}$ and bounded by $L_{X=-Y}$ in $M$. 
Set $G=G_h - \left(UM_{|\Sigma_{X=-Y}}\right)$ so that $\partial G=p(\partial(X,Y))$.
Let $\iota$ be the endomorphism of $UM$ over $M$ that maps a unit vector to the opposite one.
Set $$\begin{array}{lllll}&F(X,Y)&=[0,1/3]\times X(M) &+ \{1/3\}\times G &+ [1/3,1]\times Y(M)\\
 \mbox{and}&      F(-X,-Y)&=[0,2/3]\times (-X)(M) &+ \{2/3\}\times \iota(G) &+ [2/3,1]\times (-Y)(M).
      \end{array}
$$

Then $$F(X,Y) \cap F(-X,-Y)=[1/3,2/3] \times Y(L_{Y=-X}) - \{1/3\} \times (-X)(\Sigma_{X=-Y}) +\{2/3\} \times (Y)(\Sigma_{X=-Y}).$$
Using Lemma~\ref{lemevHM} with $F(X,X)=[0,1]\times X(M)$ to evaluate the class of $(F(X,Y) \cap F(-X,-Y))$ in $H(M)$ finishes the proof.
\eop

\noindent{\sc Proof of Proposition~\ref{proptrans}:}
Compute $lk(L_{X=Z},L_{X=-Z})$ by computing the class of $F(X,Z)\cap F(-X,-Z)$ in $H(M)$ where
$F(X,Z)$ (resp. $F(-X,-Z)$) is constructed by gluing shrinked copies of $F(X,Y)$ (resp. $F(-X,-Y)$) and $F(Y,Z)$ (resp. $F(-Y,-Z)$)
so that $[F(X,Z)\cap F(-X,-Z)]=[F(X,Y)\cap F(-X,-Y)]+[F(Y,Z)\cap F(-Y,-Z)]$.
\eop

\subsection{Proof of Theorem~\ref{thmcartorsec}}
\label{subpfcartorsec}

The first part of Theorem~\ref{thmcartorsec} follows from Lemma~\ref{lemLXYwd} and Corollary~\ref{corpontcomb}. According to Lemma~\ref{lemeultor}, two transverse sections $X$ and $Y$ are torsion sections if and only if $L_{Y=X}$ and $L_{Y=-X}$ are rationally null-homologous.
In this case, $lk(L_{Y=X},L_{Y=-X})$ only depends on the combings $[X]$ and $[Y]$ according to Lemma~\ref{lemlkind} (or to Proposition~\ref{propFXYcap}). 

Now, assume that $Y$ and $Y^{\prime}$ are such that $L_{Y=-X}$ and $L_{Y^{\prime}=-X}$ are homologous. Then $Y$ and $Y^{\prime}$ represent the same Spin$^c$-structure and there exists $k \in \ZZ$ such that $Y^{\prime}$ represents $[\gamma^{k} Y]$.
According to Corollary~\ref{coractpit}, $lk(L_{Y^{\prime}=Y},L_{Y^{\prime}=-Y})=-k$.
According to Proposition~\ref{proptrans}, 
$$lk(L_{Y^{\prime}=X},L_{Y^{\prime}=-X})-lk(L_{Y=X},L_{Y=-X})=lk(L_{Y^{\prime}=Y},L_{Y^{\prime}=-Y}).$$
Thus if 
$lk(L_{Y^{\prime}=X},L_{Y^{\prime}=-X})=lk(L_{Y=X},L_{Y=-X})$, $k=0$, and $Y$ and $Y^{\prime}$ are homotopic.
\eop

\section{On the map $p_1$}
\setcounter{equation}{0}
\label{secpone}

\subsection{The original map $p_1$ for parallelizations}
\label{subrecpone}

It has long been known that smooth compact oriented $3$-manifolds are parallelizable. 

Let $M$ be equipped with a parallelization $\tau_M\colon M \times \RR^3 \rightarrow TM$. 
Let $GL^+(\RR^3)$ denote the group of orientation-preserving linear isomorphisms of $\RR^3$.
Let $C^0((M,\partial M),(GL_+(\RR^3),\mbox{Id}))$ denote the set of maps $$g: (M, \partial M) \longrightarrow (GL^+(\RR^3),\mbox{Id})$$
from $M$ to $GL^+(\RR^3)$ that send $\partial M$ to the identity $\mbox{Id}$ of $GL^+(\RR^3)$.
Let $[(M,\partial M),(GL_+(\RR^3),\mbox{Id})]$ denote the group of homotopy classes of such maps, with the group structure induced by the multiplication of maps using the multiplication in $GL_+(\RR^3)$.

For a map $g$ in $C^0((M,\partial M),(GL_+(\RR^3),\mbox{Id}))$, define
$$\begin{array}{llll} 
\psi(g): &M \times \RR^3 &\longrightarrow  &M \times \RR^3\\
&(x,y) & \mapsto &(x,g(x)(y)).\end{array}$$
Then any parallelization $\tau$ of $M$ that coincides with $\tau_M$ on $\partial M$ reads
$$\tau = \tau_M \circ \psi(g) $$ for some $g \in C^0((M,\partial M),(GL_+(\RR^3),\mbox{Id}))$.
Thus fixing $\tau_M$ identifies the set of homotopy classes of parallelizations of $M$ fixed on $\partial M$ with
the group
$[(M,\partial M),(GL_+(\RR^3),\mbox{Id})]$.
Since $GL_+(\RR^3)$ deformation retracts onto $SO(3)$, the group $[(M,\partial M),(GL_+(\RR^3),\mbox{Id})]$ is isomorphic to $[(M,\partial M),(SO(3),\mbox{Id})]$.

The following standard proposition is proved in \cite{lesmek}.

\begin{proposition}
\label{proppone}
For any compact connected oriented $3$-manifold $M$, $[(M,\partial M),(SO(3),\mbox{Id})]$ is an abelian group, and
the degree $$\mbox{deg} \colon [(M,\partial M),(SO(3),\mbox{Id})] \longrightarrow \ZZ$$ is a group homomorphism,
that induces an isomorphism 
$$\mbox{deg} \colon [(M,\partial M),(SO(3),\mbox{Id})] \otimes_{\ZZ}\QQ \longrightarrow \QQ.$$
\end{proposition}

Let $W$ be a connected, compact $4$--dimensional manifold with signature $0$ whose boundary is
$$\partial W = \left\{
\begin{array}{ll} M \cup_{1\times \partial M} (-[0,1]\times S^2) \cup_{0 \times S^2} (-B^3) & \mbox{when}\;\partial M= S^2\\
M& \mbox{when}\;\partial M= \emptyset . \end{array}
\right. $$
When $\partial M= S^2$, $W$ has ridges and it is identified with an open subspace of one of the products $[0,1[ \times B^3$ or $]0,1] \times M$ near $\partial W$.
For any parallelization $\tau$ of $M$, the tangent vector $T_t[0,1]$ to $[0,1]$, the standard parallelization $\tau_s$ of $\RR^3$ and $\tau$ together induce a trivialization $\tau(\partial W,\tau)$ of $TW$ over $\partial W$, this trivialization reads $T_t[0,1] \oplus \tau_s$ or $T_t[0,1] \oplus \tau$.
Then the {\em Pontrjagin number\/} $p_1(\tau)$ of $\tau$ is the obstruction to extending the trivialization $\tau(\partial W,\tau) \otimes \CC$ of $TW_{|\partial W} \otimes \CC$ across $W$
(with respect to the trivialization of $det(TW)$ induced by the orientation of $W$). This obstruction lives in the module $H^4(W,\partial W;\pi_3(SU(4)))$ that is isomorphic to $\ZZ$ since $\pi_3(SU(4))=\ZZ$. For more details, see \cite[Section 1.5]{lesconst} or \cite[Proposition 6.13]{lesmek} where the following classical theorem is proved.

\begin{theorem}
\label{thmpone}
Let $M$ be a compact connected oriented $3$-manifold such that $\partial M=\emptyset\;\;\mbox{or}\;\; S^2$. For any map $g$ in $C^0((M,\partial M),(SO(3),\mbox{Id}))$, for any trivialization $\tau$ of $TM$
$$p_1( \tau \circ \psi(g))-p_1(\tau)=2\mbox{deg}(g).$$
\end{theorem}

For $n\geq 3$, a {\em spin structure\/} of a smooth $n$--manifold is a homotopy class of parallelizations over a $2$-skeleton of $M$ (that is over the complement of a point when $n=3$, if $M$ is connected).

The class of the covering map $\tilde{\rho}$ described in Notation~\ref{notrho} is the 
standard generator of $\pi_3(SO(3))=\ZZ[\tilde{\rho}]$.
The map $\rho$ can be used to describe the action of $\pi_3(SO(3))$ on the homotopy classes of parallelizations $(\tau \colon  M \times \RR^3 \rightarrow TM)$ of $M$ as follows. Let $B$ be a $3$--ball in $M$ identified with $B^3$.
Let $\tau\psi(\rho)$ coincide with $\tau$ outside $B \times \RR^3$ and read $\tau \circ \psi(\rho)$ on $B \times \RR^3$.
Set $[\tilde{\rho}][\tau]=[\tau\psi(\rho)]$. According to Theorem~\ref{thmpone}, $p_1([\tilde{\rho}][\tau])=p_1(\tau)+4$.
The set of parallelizations that induce a given spin structure form an affine space with translation group $\pi_3(SO(3))$.

The {\em Rohlin invariant\/} $\mu(M,\sigma)$ of a smooth closed $3$-manifold $M$, equipped with a spin structure $\sigma$, is the mod $16$ signature of a compact spin $4$-manifold $W$ bounded by $M$ so that the spin structure of $W$ restricts to $M$ as a stabilisation of $\sigma$. The \emph{first Betti number} of $M$ that is the dimension of $H_1(M;\QQ)$ is denoted by $\beta_1(M)$. 

Kirby and Melvin proved the following theorem \cite[Theorem 2.6]{km}.

\begin{theorem}
\label{thmim} 
For any closed oriented $3$-manifold $M$, for any parallelization $\tau$ of $M$, 
$$\left(p_1(\tau)-\mbox{dimension}(H_1(M;\ZZ/2\ZZ))-\beta_1(M)\right) \in 2 \ZZ.$$
Let $M$ be a closed $3$-manifold equipped with a given spin structure $\sigma$. Then $p_1$ is a bijection from the set of homotopy classes of parallelizations of $M$ that induce $\sigma$ to $$2\left(\mbox{dimension}(H_1(M;\ZZ/2\ZZ))+1\right) +\mu(M,\sigma)+4 \ZZ$$
When $M$ is a $\ZZ$--sphere, $p_1$ is a bijection from the set of homotopy classes of parallelizations of $M$ to $(2+4 \ZZ)$.
\end{theorem}

Extend the standard parallelization $\tau_s$ of $B^3$ as a parallelization $\widehat{\tau_s}$ of $S^3$.
When $\partial M=S^2$, form $\hat{M}=(S^3\setminus (B^3 \setminus N(\partial M)))\cup_{N(\partial M)} M$ and use $\widehat{\tau_s}$ to extend any parallelization $\tau$ of $M$ to a parallelization $\hat{\tau}$ of $\hat{M}$. 
Then it is easy to see that $p_1(\tau)=p_1(\hat{\tau})-p_1(\widehat{\tau_s})$.
In particular, according to Theorem~\ref{thmim}, $\left(p_1(\tau)-\mbox{dimension}(H_1(M;\ZZ/2\ZZ))-\beta_1(M)\right) \in 2 \ZZ$ and, when $M$ is a $\ZZ$--ball, the map $p_1$ is a bijection from the set of homotopy classes of parallelizations of $M$ to $4 \ZZ$.

\subsection{Proofs of Theorems~\ref{thmponegen} and \ref{thmimcom}}
\label{subpfmth}

\begin{lemma}
\label{lemponelk}
 Let $\tau$ be a trivialization of $TM$. Let $g \in C^0((M,\partial M),(SO(3),\mbox{Id}))$. Recall that $p_{S^2}\colon SO(3) \rightarrow S^2$ maps a transformation $t$ of $SO(3)$ to $t(N)$ where $N$ is the first basis vector of $\RR^3$. Let $X$ and $Y$ be two combings of $UM$ induced by $\tau$ and $[g][\tau]=[\tau\psi(g)]$, respectively. 
Then $$lk(L_{Y=X},L_{Y=-X})=lk((p_{S^2}\circ g)^{-1}(N),(p_{S^2}\circ g)^{-1}(-N))=-\frac12\mbox{deg}(g)$$
\end{lemma}
\bp
The first equality follows from the definition. 
It implies that $lk(L_{Y=X},L_{Y=-X})=lk((p_{S^2}\circ g)^{-1}(N),(p_{S^2}\circ g)^{-1}(-N))=lk^{\prime}(g)$ only depends on $g$. Then Proposition~\ref{proptrans} implies that $lk^{\prime}$ is a homomorphism from $[(M,\partial M),(SO(3),\mbox{Id})]$ to $\QQ$. According to Proposition~\ref{proppone} it suffices to evaluate it on the element $\rho$ viewed as a degree $2$ map of $C^0((B^3,\partial B^3),(SO(3),\mbox{Id}))$.
According to Corollary~\ref{coractpit}, when $g=\rho$, $lk(L_{Y=X},L_{Y=-X})=-1$.
\eop

\noindent{\sc Proof of Theorem~\ref{thmponegen}:}
Theorem~\ref{thmpone} and Lemma~\ref{lemponelk} show that if $X$ and $Y$ extend to parallelizations $\tau(X)$ and $\tau(Y)$, then
$$p_1(\tau(Y))-p_1(\tau(X))=-4lk(L_{Y=X},L_{Y=-X}).$$
For any torsion combing $[Y]$, define $p_1([Y])$ from a combing $[X]$ that extends to a parallelization by
$$p_1([Y])=p_1([X])+4lk(L_{X=Y},L_{X=-Y}).$$
Thanks to Proposition~\ref{proptrans}, since this formula is valid for combings that extend to parallelizations, this definition does not depend on the choice of $X$. Now, Proposition~\ref{proptrans} implies that the above formula is valid for all pairs of torsion combings.

Since $[-X]=[X]$ for a section $X$ that extends as a trivialization, we deduce that $p_1([-Y])=p_1([Y])$, for all torsion sections $Y$, from the above definition.

According to the following Lemma~\ref{lemactpithree}, Proposition~\ref{propSpinstruc} ensures the injectivity of the restriction of $p_1$ to any torsion Spin$^c$-structure.
\eop

\begin{lemma}
 \label{lemactpithree}
For any torsion combing $[X]$, $p_1(\gamma[X])-p_1([X])=4$.
\end{lemma}
Recall Corollary~\ref{coractpit}. \eop

\begin{proposition}
\label{propthmimcom}
With the notation of Theorem~\ref{thmpont}, if $(L,L_{\parallel})$ is a framed rationally null-homologous link of the interior of $M$, then
$$p_1(C(\tau,L,L_{\parallel}))=p_1(\tau)-4 lk(L,L_{\parallel}).$$
\end{proposition}
\bp Assume that $\tau$ reads $(X,X_2,X_3)$ so that $L=L_{Y=-X}$. Then $X$ and $X_2$ are homotopic sections of $UM$
so that $p_1(\tau)=p_1(X)=p_1(X_2)$ and, according to Theorem~\ref{thmponegen}, $p_1(Y)=p_1(\tau)-4 lk(L_{Y=X_2},L_{Y=-X_2})$. The link $(L_{Y=X_2},L_{Y=-X_2})$ is isotopic to  $(L,L_{\parallel})$.
\eop

\noindent {\sc Proof of Theorem~\ref{thmimcom}:}
According to Theorem~\ref{thmpont}, any torsion combing $Y$ is homotopic to $C(\tau,L,L_{\parallel})$, for some $L$ and $\tau$ as in Proposition~\ref{propthmimcom} and in its proof. In particular, since $L=L_{Y=-X}$, $p_1(Y)\in p_1(\tau)-4 \ell([L_{Y=-X}]))$ and $p_1(Y)\in  p_1(\tau) -4 \ell(\mbox{\rm Torsion}(H_1(M;\ZZ))$.
Conversely, any element in $\ell(\mbox{\rm Torsion}(H_1(M;\ZZ))$ reads $lk(L,L_{\parallel})$ for some rationally null-homologous link link $L$.
\eop

\subsection{Identifying $p_1$ with the Gompf invariant}
\label{secGo}

Let us first recall the definition of the Gompf invariant.
An {\em almost-complex structure\/} on a smooth $4$-dimensional manifold $W$ is an operator $J$ such that $J^2=-\mbox{Id}$, acting smoothly on the tangent space to $W$, fiberwise.
An almost-complex structure on $W$ induces a combing of $\partial W$, that is the class of the image $[JN=J(N(\partial W))]$ under $J$ of the outward normal $N(\partial W)$ to $W$. Gompf showed that all the combings of a $3$-manifold appear as combings $JN$ for some $W$ \cite[Lemma~4.4]{Go}, this will be reproved below. The {\em first Chern class\/} $c_1(TW,J)$ of $(TW,J)$ is the obstruction to trivializing $TW$ over the two-skeleton of $W$ as an almost-complex manifold (the induced trivialization of $TW$ must read $(X,JX,Y,JY)$). The class $c_1(TW,J)$ lives in $H^2(W;\ZZ)$. (The first Chern class $c_1$ of a complex vector bundle is the Euler class of the corresponding determinant bundle. The reader can check that the definitions coincide in this case.)  Its restriction to $H^2(\partial W;\ZZ)$ is $e(JN^{\perp})$ so that
the boundary of the Poincar\'e dual $Pc_1(TW,J)$ of $c_1(TW,J)$ is Poincar\'e dual to $e(JN^{\perp})$. When $JN$ is a torsion combing, this boundary $\partial Pc_1(TW,J)$ is a torsion element of $H_1(\partial W;\ZZ)$ so that there exists a rational $2$--chain $\Sigma$ of $\partial W$ such that $(Pc_1(TW,J) \cup \Sigma)$ is a closed rational $2$-cycle of $W$. The algebraic self-intersection of this rational cycle is independent of $\Sigma$ and it is denoted by $(Pc_1(TW,J))^2$, and the \emph{Gompf invariant} $\theta_G(JN)$ that is denoted by $\theta(JN)$ in \cite[Section 4]{Go}
is $$\theta_G(JN)=(Pc_1(TW,J))^2-2\chi(W)-3\,\mbox{signature}(W)$$
where $\chi$ stands for the Euler characteristic.

In this subsection, we prove that $\theta_G=p_1$.

\begin{lemma}
\label{lemeqfram}
When a combing $X$ of $M$ extends as a parallelization, $\theta_G([X])=p_1([X])$.
\end{lemma}
\bp For a rank $2k$ complex bundle $\omega$ seen as a rank $4k$ real bundle $\omega_{\RR}$, $p_1(\omega_{\RR})=c_1^2(\omega)-2c_2(\omega)$, where $c_2$ denotes the second Chern class that is the Euler class of $\omega_{\RR}$ for a rank $2$ complex bundle $\omega$. See \cite[Definition p.158 \S~14 and Corollary~15.5]{milnorsta}.
Let $(W,J)$ be an almost-complex connected compact manifold bounded by $M$ such that $X=JN$, let $Y$ be a nowhere zero section of $X^{\perp} \subset TM$.
Consider the almost-complex parallelization $(N,Y)$ inducing the real parallelization $(N,JN,Y,JY)$ of $TW_{|M}$, and the complex bundle $\omega$ over $(W \cup_M (- W))$
that is trivial with fiber $\CC N \oplus \CC Y$ over $(-W)$ and that coincides with the initial one over $W$.
Since the characteristic classes $p_1$, $c_1$ and $c_2$ of $\omega_{\RR}$ or $\omega$ trivially restrict to $H^{\ast}(-W)$, they come from classes of $H^{\ast}(W \cup_M (- W),-W) \cong H^{\ast}(W,M)$. Thus $p_1(\omega_{\RR})$ is the image of $p_1(W,(JN,Y,JY))[W,\partial W] \in H^4(W,\partial W)$, and $c_2(\omega)$ is the image of $c_2(TW,N) \in H^4(W,\partial W)$ that is 
$\chi(W)[W,\partial W]$ since $c_2$ is the obstruction to extending $N$ as a nowhere zero section of $TW$, that is the relative Euler class of $(TW,N)$.
Similarly, $c_1(\omega)$ is the image of a lift $\tilde{c}_1$ of $c_1(TW,J)$ in $H^2(W,\partial W)$, where $P\tilde{c}_1$ is represented by a cycle of $W$ that can be constructed as in the definition of $(Pc_1(TW,J))^2$ before Lemma~\ref{lemeqfram}. The Poincar\'e dual $Pc_1(\omega)$ of $c_1(\omega)$ is the image of this cycle in $H_2(W \cup_M (- W))$ and
$p_1(W,(JN,Y,JY))=(Pc_1(TW,J))^2-2\chi(W)$.
\eop

\begin{lemma}
\label{lemeqcomfram}
When a combing $X$ of $M$ extends as a parallelization, $\theta_G([\gamma X])=\theta_G([X])+4$.
\end{lemma}
\bp According to Lemma~\ref{lemactpithree}, $p_1([\gamma X])=p_1([X])+4$ for any $[X]$. \eop

Any closed oriented connected $3$-manifold $M$ is the boundary of a $4$-manifold $$W_L=B^4 \bigcup_{L \times D^2 \subset S^3}\coprod_{i=1,\dots,n}(D^2 \times D^2)^{(i)}$$
obtained from $B^4$ by attaching $2$-handles $(D^2 \times D^2)^{(i)}_{i=1,\dots,n}$ along a tubular neighborhood $L \times D^2$ of a framed link $L=(K_i,\mu_i)_{i=1,\dots,n}$.
Such a framed link $L$ is an {\em integral surgery presentation\/} of $W_L$ and $M$. The $K_i$ are the components of $L$, the $\mu_i$ are the surgery parallels $K_i \times \{1\} \subset K_i \times D^2$ that frame the $K_i$, and the handle $(D^2 \times D^2)^{(i)}$ is attached by a natural identification of $K_i \times D^2 \subset \partial B^4$ with $((-S^1) \times D^2)^{(i)}$
that restricts to $\mu_i$ as an orientation-reversing homeomorphism onto $(S^1 \times \{1\})^{(i)}$.

According to Kaplan \cite{Ka}, we can furthermore demand that $lk(K_i,\mu_i)$ is even for any $i$, in the statement above. In this case, we shall say that the surgery presentation is {\em even.\/}

\begin{lemma}
Let $L$ be an even surgery presentation of $M$. 
There is an almost-complex structure $J^0$ on $W_L$ (described below) such that $e(J^0N^{\perp})=0$.
For any $Spin^c$ structure $\xi$ on $M$, there is at least one almost complex structure $J$ on $W_L$ (described below) such that the class of $JN$ belongs to $\xi$ and, if $JN$ is a torsion combing, then $p_1(JN)-p_1(J^0N)=\theta_G(JN)-\theta_G(J^0N)$.
\end{lemma}
\bp
We shall only consider almost-complex structures $J$ that are {\em compatible with a\/} given {\em Riemannian metric\/} in the following sense:
$J$ preserves the Riemannian metric and $Jx$ is orthogonal to $x$ for any $x$.
Our almost-complex structures $J$ of $4$-manifolds also induce the orientation via local parallelizations of the form $(X,JX,Y,JY)$.
Below, $B^4$ is seen as the unit complex ball of $\CC^2$, it is equipped with its usual Riemannian structure.
It is also seen as the unit ball of the quaternion field $\HH=\CC \oplus \CC j$, so that $S^3$ is identified with the group of unit quaternions and $T_xS^3$ is the space of quaternions orthogonal to $x$.

A homotopy $$\begin{array}{llll}JN \colon&[-1,0]\times S^3 &\rightarrow &TS^3 \\&(t,x) &\mapsto &JN(t,x) \in T_xS^3\end{array}$$ such that $JN(-1,x)=ix$, and $\left\Vert JN(t,x)\right\Vert=1$  induces a homotopic almost-complex structure on $B^4$ as follows, the complex structure is unchanged outside a collar $[-1,0]\times S^3$ of the boundary of $B^4$, and the operator $J$ of the almost-complex structure maps the unit tangent vector to $[0,1] \times \{x \in S^3\}$ at $(t,x)$ to $JN(t,x)$. Note that $J$ is completely determined by these conditions.
If such a homotopy is such that $JN(0,.)$ is tangent to $K_i \times \{y\}$ on $K_i \times D^2$, then the associated almost-complex structure
$J$ preserves the tangent space to $\{x\} \times D^2$ and it uniquely extends to $(D^2 \times D^2)^{(i)}$ so that $J$ preserves the tangent space to $\{x\} \times D^2$ and $J$ is compatible with the product Riemannian structure on $(D^2 \times D^2)^{(i)}$.
In particular $J$ maps the outward normal to $(D^2 \times S^1)^{(i)} \subset M$ at $(x,y \in S^1)$ to the unit tangent vector to $(\{x\} \times S^1)^{(i)}$ at $(x,y)$.

Before smoothing the ridges, $W_L$ reads $(\RR^2 \setminus \{(x,y);x<-1,y>-1\} ) \times (-K_i) \times S^1$ near $K_i \times S^1$.
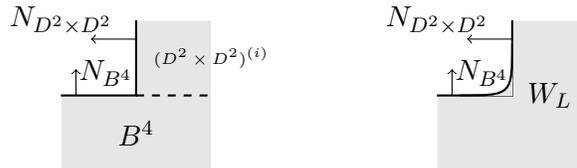
\begin{figure}[h]
\begin{center}
\begin{tikzpicture}
\useasboundingbox (-1.1,-1.1) rectangle (6.1,1.1);
\draw [fill=gray!20,draw=white] (-1,-1) rectangle (1,1);
\draw [fill=white,draw=white] (-1,0) rectangle (0,1);
\draw [thick] (-1,0) -- (0,0) -- (0,1);
\draw [thick,dashed] (0,0) -- (1,0);
\draw [->] (-.8,0)  -- (-.8,.35);
\draw [->] (0,.75) -- (-.6,.75);
\draw (-.4,.3) node{\small $N_{B^4}$} (-1,1) node{\small $N_{D^2 \times D^2}$} (0,-.5) node{\small $B^4$} (1,.5) node{\tiny $(D^2 \times D^2)^{(i)}$};
\draw [fill=gray!20,draw=white] (4,-1) rectangle (6,1);
\draw [fill=white,draw=white] (4,0) rectangle (5,1);
\draw [thick] (4,0) -- (5,0) -- (5,1);
\draw [fill=gray!20,draw=white] (4.3,0) .. controls (5,0) .. (5,.7) -- (5,0) -- (4.3,0);
\draw [thick] (4.3,0) .. controls (5,0) .. (5,.7);
\draw [->] (4.2,0)  -- (4.2,.35);
\draw [->] (5,.75) -- (4.4,.75);
\draw (4.6,.3) node{\small $N_{B^4}$} (4,1) node{\small $N_{D^2 \times D^2}$} (5.5,0) node{\small $W_L$};
\end{tikzpicture}
\caption{$W_L$ near $K_i \times S^1$ before and after smoothing.}
\label{figwl}
\end{center}
\end{figure}
The $4$-manifold $W_L$ is next smoothed around $ K_i \times S^1$, the smoothing adds the product of $ K_i \times S^1$ by a triangle with two orthogonal straight sides and a smooth hypothenuse that makes null angles with the two straight sides.  See Figure~\ref{figwl}.

This new piece may be seen as a part of a $D^2 \times \RR^2$ that contains $D^2 \times D^2$,
so that $J$ naturally extends there.

In the plane of the triangle, the normal $N$ reads $N=\cos(\theta)N_{B^4} + \sin(\theta)N_{D^2 \times D^2}$ for some $\theta \in [0,\pi/2]$, so that
$JN$ reads $JN=\cos(\theta)JN_{B^4} + \sin(\theta)JN_{D^2 \times D^2}$ and $JN$ goes from the tangent to $K_i \times \{y\}$ to the tangent to $(\{x\} \times S^1)^{(i)}$ on $T_{(x,y)}K_i \times S^1$ by the shortest possible way on the smooth hypothenuse.

Then $J$ and $JN$ are completely determined on $W_L$ by the homotopy $JN$ on $[-1,0]\times S^3$, and we now study them as a function of this homotopy.

We shall consider homotopies induced by homotopies of orthonormal parallelizations, i.e. homotopies $JN$ such that there is a homotopy
$V \colon[-1,0]\times S^3 \rightarrow T_xS^3$ where $V(t,x) \in T_xS^3$, $V(t,x) \perp JN(t,x)$, $\left\Vert V(t,x) \right\Vert=1$ and $V(-1,x)=jx$.
Furthermore, our homotopies are such that $JN(0,.)$ is tangent to $K_i \times \{y\}$ on $K_i \times D^2$, so that
$V(0,x)$ induces a framing of $K_i$. The linking number of $K_i$ with the parallel of $K_i$ induced by this framing is denoted by $r_i$. Recall that $H_1(SO(3);\ZZ)=\pi_1(SO(3))=\ZZ/2\ZZ$ is generated by a loop of rotation $(\exp(i\theta) \mapsto \rho(\theta,A))$ for a fixed arbitrary axis $A$. 

\begin{sublemma}
\label{sublriodd}
The integers $r_i$ are odd.
\end{sublemma}
\noindent{\sc Proof of Sublemma~\ref{sublriodd}: }
Let $\Sigma$ be a Seifert surface of $K_i$, then $TM_{|\Sigma}$ has a trivialization $\tau_{\Sigma}$ whose third vector is the positive normal $N\Sigma$ to $\Sigma$, and whose first vector over $K_i$ is obtained from the tangent vector $v_K$ to $K_i$ by rotating it $(-\chi(\Sigma))$ times around the axis $N\Sigma$, along $K_i$. On the other hand, the first vector of the restriction to $K_i$ of the trivialization $\tau_{JV}$ induced by $JN(0,.)$ and $V(0,.)$ is $v_K$ and its third vector is obtained from $N\Sigma$ by rotating it $r_i$ times around $v_K$ along $K_i$. Then $\tau_{\Sigma}^{-1}\circ \tau_{JV}$ induces a map from $\Sigma$ to $SO(3)$ whose restriction to $K_i$ represents a trivial homology class in $H_1(SO(3))$. Since the class of this restriction is $(r_i+\chi(\Sigma))$ mod $2$ and since $\chi(\Sigma)$ is odd, $r_i$ is odd, too.
\eop

\begin{sublemma}
\label{sublrioddarb}
The integers $r_i$ may be changed to any arbitrary odd number, by perturbing the homotopy near $K_i \times D^2$.
\end{sublemma}
\noindent{\sc Proof of Sublemma~\ref{sublrioddarb}: }
Assume without loss that $JN(0,.)$ is tangent to $K_i \times \{y\}$ on a bigger tubular neighborhood $K_i \times 2D^2$. Let $e_1$ denote the first basis vector of $\RR^3$. Consider a map $$\begin{array}{lllll}F\colon &[0,1]\times \frac{\RR}{2\pi\ZZ} &\rightarrow &SO(3) &\\
& (t,\theta) & \mapsto &  \mbox{Id}   &\;\mbox{if}\; t=1    \; \mbox{or}\;  \theta \in 2\pi\ZZ    \\                                                                                                                           
&&&  \rho(2\theta,e_1) & \; \mbox{if}\; t=0.                                                                                                                                   \end{array}
$$ Then $(JN,V,JV)(0,.)$ may be replaced on $K_i \times 2D^2$, by the homotopic 
$$(0,(\exp(i\theta), u \exp(i\eta))) \mapsto 
F(\max(0,u-1),k_i\theta)\left((JN,V,JV)(0,(\exp(i\theta), u \exp(i\eta)))\right)$$ for some integer $k_i$. Since this changes $r_i$ to $r_i+2k_i$, this shows that $r_i$ can be changed to any odd number.
\eop

Now, the obstruction to extending $V$ as a unit vector tangent to the second almost-complex factor $D^2$ across $(D^2 \times .)^{(i)}$ is $-(r_i-lk(\mu_i,K_i))$,
and the obstruction to extending $JN$ that is the tangent to $K_i \times \{y\}$ as a unit vector tangent to the first almost-complex factor across $(D^2 \times .)^{(i)}$ is $1$. In particular, the Poincar\'e dual of the Chern class $c_1(TW_L,J)$ may be represented by a chain that does not intersect $B^4$ and that intersects $(D^2 \times D^2)^{(i)}$ as $(1-r_i+lk(\mu_i,K_i))(0\times D^2)^{(i)}$.
Let $J^0N$ be a homotopy such that $r_i=lk(\mu_i,K_i)+ 1$. Then $c_1(TW_L,J^0)=0$ since $H^2(B^4,S^3)=0$ and $\theta_G(J^0N)=-2\chi(W_L)-3\;\mbox{signature}(W_L).$

Change $r_i$ to $(r_i+2k_i)$ as in Sublemma~\ref{sublrioddarb} and denote
the obtained almost-complex structure by $J$. Compare the induced vector fields and compute $L_{JN =V^0}$ and $L_{JN =-V^0}$.

Note that composing the map $F$ in the proof of Sublemma~\ref{sublrioddarb} by the evaluation $p_{S^2}$ at $e_1$ provides a degree $\pm 1$ map from $([0,1]\times S^1,\partial [0,1]\times S^1)$ to $(S^2,e_1)$.
Thus there exists a well-determined $\varepsilon =\pm 1$ such that  $L_{JN =V^0}$ and $L_{JN =-V^0}$ are homologous to $\varepsilon \sum_{i=1}^n k_i m_i$ in $(L \times D^2) \setminus L$ where $m_i$ is a meridian of $K_i$. We can furthermore assume that $L_{JN =V^0} \subset L \times u S^1$ and $L_{JN =-V^0} \subset L \times u^{\prime} S^1$ for two distinct elements $u$ and $u^{\prime}$ of $]1,2[$. Let $m_i$ (resp. $m_{i\parallel}$) denote a meridian of $K_i$ in $L \times u S^1$ (resp. in $L \times u^{\prime} S^1$).

Since the meridians $m_i$ generate $H_1(M;\ZZ)$, for any $Spin^c$-structure $\xi$, there exists an almost-complex structure $J$ as above such that $JN$ belongs to $\xi$.
The combing $JN$ is torsion if and only if $L_{JN =-V^0}$ represents a torsion element in $H_1(M;\ZZ)$. Assume that $JN$ is torsion from now on. According to Theorem~\ref{thmponegen}
$$p_1(JN)-p_1(J^0N) =p_1(JN)-p_1(V^0)=-4lk(\sum_{i=1}^n k_i m_i,\sum_{i=1}^n k_i m_{i\parallel}).$$
On the other hand, since the boundary of 
$Pc_1(TW_L,J)$ is homologous to $2 L_{JN =-V^0}$,
$Pc_1(TW_L,J)$ is represented by $2 \varepsilon\sum_{i=1}^n k_i (0\times D^2)^{(i)}$.
In order to compute $(Pc_1(TW_L,J))^2$, consider a parallel copy $Pc_1(TW_L,J)_{\parallel}=2 \varepsilon\sum_{i=1}^n k_i (x\times D^2)^{(i)}$, and let $(-\partial Pc_1(TW_L,J))$ and $(-\partial Pc_1(TW_L,J)_{\parallel})$ bound $\Sigma$ and $\Sigma_{\parallel}$ in $M$, respectively, so that
$$\begin{array}{ll}\theta_G(JN)-\theta_G(J^0N)&=(Pc_1(TW,J))^2\\&=\langle 2 \varepsilon\sum_{i=1}^n k_i(0\times D^2)^{(i)} \cup \Sigma ,2 \varepsilon\sum_{i=1}^n k_i(x\times D^2)^{(i)} \cup \Sigma_{\parallel}\rangle_{W_L \cup_{\partial W_L = 0 \times M} [0,1] \times M}
\\&=\langle (-[0,1/2] \times \partial \Sigma) \cup (1/2 \times \Sigma), (-[0,2/3] \times \partial \Sigma_{\parallel}) \cup (2/3 \times \Sigma_{\parallel}) \rangle_{[0,1]\times M}\\&=-\langle \Sigma,  \partial \Sigma_{\parallel} \rangle_{M}\\&=p_1(JN)-p_1(J^0N).\end{array}$$
\eop

The previous lemma, Lemma~\ref{lemactpithree} and the transitivity of the action of $\pi_3(S^2)$ on the combings of a $Spin^c$-structure reduce the proof that $\theta_G=p_1$ to the proof of the following lemma.

\begin{lemma}
 $\theta_G([\gamma X])-\theta_G([X])=4$ for any combing $[X]$.
\end{lemma}
\bp
We refer to the previous proof.
Add a trivial knot $U$ framed by $+1$ to a surgery presentation $L$, such that $W_L$ is equipped with an almost-complex structure $J$.
The structure $J$ is homotopic to a structure $J^{(1)}$ that extends on $W_{L\cup U}$ so that $Pc_1(TW,J^{(1)})$ is $ (0\times D^2)^{(0)}$. Then $\theta_G(J^{(1)}N)-\theta_G(JN)=1-2-3=-4$. The structure $J$ is also homotopic to a structure $J^{(3)}$ that extends on $W_{L\cup U}$ so that $Pc_1(TW,J^{(3)})$ is $3 (0\times D^2)^{(0)}$, then $\theta_G(J^{(3)}N)-\theta_G(J^{(0)}N)=9-2-3=4$.
These two combing modifications sit in a $3$-ball of $M$, so that each of them corresponds to the action of an element of $\pi_3(S^2)$ independent of $(M,J)$.
According to Lemma~\ref{lemeqcomfram}, $[J^{(1)}N]=[\gamma^{-1}JN]$ and $[J^{(3)}N]=[\gamma JN]$.
Since the above process allows us to inductively represent all the combings $[\gamma^k JN]$, by adding some disjoint trivial knots framed by $+1$, and to prove that $\theta_G(\gamma^k JN)-\theta_G(\gamma^{k-1} JN)=4$, for all $k \in \ZZ$, we are done. \eop

\begin{remark}
For a natural integer $k$ and for a surgery presentation $L$ of $M$ in $S^3$, let $L(k)$ be the surgery presentation of $M$ obtained from $L$ by adding $k$ trivial knots framed by $+1$.
 On our way, we have proved that for any combing $[X]$, for any even surgery presentation $L$ of $M$, there exists a natural integer $k$ and an almost complex structure $J$ on $W_{L(k)}$ such that $[X]=[JN]$.
\end{remark}

\subsection{More variations of $p_1$}
\label{subonemore}

In applications, combing modifications often arise as in Definition~\ref{defchangecomb} or as in the statement of Proposition~\ref{propLcombbis} below. We show how the variation formula of Theorem~\ref{thmponegen} applies in these settings to yield other useful variation formulas.

\begin{lemma}
 \label{lemEul}
Let $M$ be equipped with a torsion combing $X$.
Let $L$ be a rationally null-homologous link in the interior of $M$. Let $Z$ be a section of\/ $UM$ orthogonal to $X$, such that $Z$ is defined on $L$ and $\partial M$.
Extend $Z$ as a section $\tilde{Z}$ of the $D^2$-bundle $X^{\perp}$, so that $\tilde{Z}$ is transverse to the zero section.
Let $L(Z \subset X^{\perp})$ be the zero locus of $\tilde{Z}$ cooriented by the fiber $D^2$ of $X^{\perp}$.
Then $L(Z \subset X^{\perp})$ is a link of $M \setminus L$ that represents the Poincar\'e dual of the relative Euler class of $(X^{\perp},Z)$, and $L(Z \subset X^{\perp})$ is homologous to the Poincar\'e dual of $e(X^{\perp})$.
\end{lemma}
\eop
\begin{remark}
Lemma~\ref{lemEul} can be taken as a definition of the relative Euler class in this case.
The obstruction to extending $Z$ across a $2$--cycle of $(M,L \cup \partial M)$ is the intersection of the $2$--cycle with $L(Z \subset X^{\perp})$.
\end{remark}

\begin{definition}
\label{defchangecomb}
Let $X$ be a section of $UM$.
Let $L$ be a link in the interior of $M$ and let $Z$ be a section of $UM_{|L}$ orthogonal to $X$. Let $\eta=\pm 1$, let $L_{\parallel}$ be a parallel of $L$ and let $N(L)$ be a tubular neighborhood of $L$ where $Z$ is extended as a section of $UM$ orthogonal to $X$.
Let $\rho(\theta,X)$ denote the rotation with axis $X$ and angle $\theta$. Let $D^2=\{u \exp(i\theta); u\in [0,1], \theta \in [0,2\pi]\}$ be the unit disk of $\CC$.
Define $D(X,L,L_{\parallel},Z,\eta)$ (up to homotopy) as the section of $UM$ that coincides with $X$ outside $N(L)$ and that reads as follows in $N(L)$ that is trivialized with respect to $L_{\parallel}$ so that it reads $D^2\times L$.
\begin{itemize}
\item $D(X,L,L_{\parallel},Z,\eta)(0,k \in L)=-X(0,k)$,
\item when $u \in]0,1]$, $[-X, D(X,L,L_{\parallel},Z)( u \exp(i\theta),k)]$ is the geodesic arc of length $u\pi$ of the half great circle $[-X,X]_{\rho(\eta\theta,X)(Z)}$ from $(-X)$ to $X$ through $\rho(\eta\theta,X)(Z)$, where $X$ and $Z$ stand for $X(u \exp(i\theta),k)$ and $Z(u \exp(i\theta),k)$, respectively,
\end{itemize} so that $D(X,L,L_{\parallel},Z,\eta)(1/2,k)=Z(1/2,k)$.
Note that the homotopy class of $D(X,L,L_{\parallel},Z,\eta)$ can also be defined by the following formula. $$D(X,L,L_{\parallel},Z,\eta)(u \exp(i\theta),k)=\rho(\pi(1+u),\rho(\eta\theta-\pi/2,X)(Z))(X)(u \exp(i\theta),k).$$
\end{definition}

\begin{remark}
\label{rkfram}
Note that with the notation of Remark~\ref{rkxfram} $$C(X,L,\sigma)=D(X,L,L_{\parallel},Z(\sigma,\sigma_N),-1).$$
\end{remark}

\begin{proposition}
\label{propLcomb}
Under the assumptions of Lemma~\ref{lemEul} above, let $\eta=\pm 1$, let $L_{\parallel}$ be a parallel of $L$ and let $N(L)$ be a tubular neighborhood of $L$ where $Z$ is extended as a section of $UM$ orthogonal to $X$. For the combing $D(X,L,L_{\parallel},Z,\eta)$ of Definition~\ref{defchangecomb}, $$p_1(D(X,L,L_{\parallel},Z,\eta)) - p_1(X) = 4 lk(L, \eta L(Z \subset X^{\perp})- L_{\parallel}).$$
\end{proposition}
\bp
Set $Y=D(X,L,L_{\parallel},Z,\eta)$. Let $\tau$
be the parallelization of $N(L)$ with first vector $X$ and second vector $Z$. Then $\tau^{-1}$ maps $Y( D^2/\partial D^2 \times k)$ to the sphere $S^2$ with degree $(-\eta)$
so that $L_{Y=-X}=-\eta L$ and $L_{X=-Y}=\eta L$.
In order to use Theorem~\ref{thmponegen}, deform $X$ to $\tilde{X}$ to make it transverse to $Y$ using $\tilde{Z}$ as follows.
Let $N_{1/3}(L)=\{(u \exp(i\theta), k \in L) \in N(L);u \in [0,1/3]\}$ and $N_{2/3}(L)=\{(u \exp(i\theta),k) \in N(L);u \in [0,2/3]\}$.
Consider a function $\chi \colon M \rightarrow [0,1]$ that maps $\left(M \setminus N_{2/3}(L)\right)$ to $1$ and $N_{1/3}(L)$ to $0$.
Let $\varepsilon$ be a very small positive real number, set $\tilde{X}=\frac1{\parallel X + \varepsilon\chi\tilde{Z}\parallel}(X + \varepsilon\chi\tilde{Z})$ so that $\tilde{X}(M)$ is now transverse to 
$Y(M)$.
Outside $UM_{|N(L)}$, $\tilde{X}(M) \cap Y(M)$ reads $Y(L(Z \subset X^{\perp}))$, whereas
on $UM_{|N(L)}$, $Y(M) \cap \tilde{X}(M)$  reads $Y(-\eta L_{\parallel})$ because $Y$ covers $S^2$ with degree $(-\eta)$ along a fiber of $N(L)$.
\eop

We have the two immediate corollaries.
\begin{corollary}
\label{corone}
Under the hypotheses of Proposition~\ref{propLcomb}, 
when $Z$ extends as a section of the unit bundle of $X^{\perp}$ on $M$,
$$p_1(D(X,L,L_{\parallel},Z,\eta)) = p_1(X)- 4 lk(L,L_{\parallel}).$$\end{corollary}

\begin{corollary}
Under the hypotheses of Proposition~\ref{propLcomb}, let $K=\{K(\exp(i\kappa) \in S^1)\}$ be a component of $L$, let $r \in \ZZ$, and let $Z_r=Z$ on $L \setminus K$ and $Z_r(K(k=\exp(i\kappa)))=\rho(r\kappa,X)(Z)(k)$.
Then $$p_1(D(X,L,L_{\parallel},Z_r,\eta)) - p_1(D(X,L,L_{\parallel},Z,\eta))=4 \eta r.$$
\end{corollary}

Note that under the hypotheses of Proposition~\ref{propLcomb}, when $X$ is tangent to $L$, if $Z$ is induced by $L_{\parallel}$, then $D(X,L,L_{\parallel},Z,1)$ is independent of $Z$ and $L_{\parallel}$.

\medskip

The following combing modification also arises in the study of combings associated with  Heegaard diagrams.

\begin{proposition}
\label{propLcombbis}
Let $M$ be equipped with a torsion combing $X$.
Let $N_0(L)$ denote a tubular neighborhood of a rationally null-homologous link $L$ in the interior of $M$. Let $L_{2} \subset \partial N_0(L)$ be a satellite of $L$ such that the restriction to $L_2$ of the bundle projection of $N_0(L)$ onto $L$ defines a $2$-fold covering of $L$. Let $s$ be the involution of $L_2$ that exchanges two points in a fiber of this covering. 
Pick a parallelization $\tau$ of $M$ such that $X$ is constant with respect to $\tau$ over $N(L)$. Let $Z$ be a section orthogonal to $X$ of the restriction of $UM$ to $L_2$, such that $Z(s(k))=-Z(k)$. Define $D(X,L,L_2,Z,-1)$ as follows.
On intervals $I$ of $L$, trivialize a larger tubular neighborhood $N(L)$ ($N_0(L) \subset N(L)$) as $D^2 \times I$ so that $(D^2 \times I) \cap L_2$
reads $\{-1/2,1/2\} \times I$, and define $D(X,L,L_2,Z,-1)$ as in Proposition~\ref{propLcomb} on these portions:
\begin{itemize}
\item $D(X,L,L_2,Z,-1)(0,k \in I)=-X(0,k)$,
\item when $u \in]0,1]$, $[-X, D(X,L,L_2,Z)( u \exp(i\theta),k)]$ is the geodesic arc of length $u\pi$ of the half great circle $[-X,X]_{\rho(-\theta,X)(Z(1/2,k))}$ from $(-X)$ to $X$ through $\rho(-\theta,X)(Z(1/2,k))$,
\end{itemize} so that $D(X,L,L_2,Z,-1)(1/2,k)=Z(1/2,k)$.
Let $f$ be a smooth increasing surjective function from an interval $I$ to $[0,\pi]$, such that all derivatives of $f$ vanish at the ends of $I$. Let $k \in \ZZ$.
Define $$\begin{array}{llll}T^k \colon & D^2 \times I & \longrightarrow & D^2 \times I\\
      &(u\exp(i\theta),t) & \mapsto & (u\exp(i(\theta+kf(t))),t)
     \end{array}$$ so that $T$ is a half-twist.
Assume that $D^2 \times I$ is a part of $N(L)$ as above and let
$(T^k(L_2),T^k_{\ast}(Z))$ coincide with $(L_2,Z)$ outside $D^2 \times I$ and read $(T^k(L_2),T^k_{\ast}(Z))$ on  $D^2 \times I$ where, for $\theta \in \{-1/2,1/2\}$, $$T^k_{\ast}(Z)((\exp(i(\theta+kf(t)))/2,t))=\rho(kf(t),X)(Z((\exp(i\theta)/2,t))).$$
Then $$p_1(D(X,L,T^k(L_2),T^k_{\ast}(Z),-1))-p_1(D(X,L,L_2,Z,-1))=-4k $$
\end{proposition}
\bp
The variation of a combing under some $T^k$ sits inside the ball $D^2 \times I$. Therefore the corresponding variation of $p_1$ may be read in this ball. It does not depend on the trivialization of the ball induced by $X$ and $Z$, since all of them are homotopic. Therefore, it only depends on $k$, linearly. The coefficient is obtained by looking at the effect of the twist $T^2$ on a $D(X,L,L_{\parallel},Z,-1)$ as in Proposition~\ref{propLcomb}.
\eop

\section{The $\Theta$-invariant of combings}
\setcounter{equation}{0}
\label{sectheta}

In this section, we present a self-contained homogeneous definition of an invariant $\Theta$ of combings of rational homology balls. This definition is deeply inspired from the definition of $\Theta_{KKT}$ that can be found in \cite[Section 6.5]{lessumgen} and in \cite[Theorem 2.14]{lesbonn}.

\subsection{On configuration spaces}

Recall that {\em blowing up\/} a submanifold $A$ means replacing it by its unit normal bundle. See Definition~\ref{defblowup}.

In a closed $3$-manifold $R$, we fix a point $\infty$ and define
$C_1(R)$ as the compact $3$-manifold obtained from $R$ by blowing up $\{\infty\}$. This space $C_1(R)$ is a compactification of $\check{R}=(R \setminus \{\infty\})$.

The {\em configuration space\/} $C_2(R)$ is the compact $6$--manifold with boundary and corners obtained from $R^2$ by blowing up $(\infty,\infty)$, and the closures of $\{\infty\} \times \check{R}$, $\check{R} \times \{\infty\}$ and the diagonal of $\check{R}^2$, successively.

Then $\partial C_2(R)$ contains the unit normal bundle  to the diagonal of $\check{R}^2$. This bundle is canonically isomorphic to $U\check{R}$ via the map 
$$[(x,y)] \in \frac{\frac{T_r\check{R}^2}{\mbox{\tiny diag}} \setminus \{0\}}{\RR^{\ast +}} \mapsto [y-x] \in \frac{T_r\check{R} \setminus \{0\}}{\RR^{\ast +}}.$$

Since $((\RR^3)^2 \setminus \mbox{diag})$ is homeomorphic to $\RR^3 \times ]0,\infty[ \times S^2$ via the map $$(x,y) \mapsto (x,\parallel y-x \parallel, \frac{1}{\parallel y-x \parallel}(y-x)),$$ $((\RR^3)^2 \setminus \mbox{diag})$ is homotopy equivalent to $S^2$. In general, $C_2(R)$ is homotopy equivalent to $(\check{R}^2 \setminus \mbox{diag})$. When $R$ is a rational homology sphere, $\check{R}$ is a rational homology $\RR^3$ and the rational homology of $(\check{R}^2 \setminus \mbox{diag})$ is isomorphic to the rational homology of $((\RR^3)^2 \setminus \mbox{diag})$. Thus, $C_2(R)$ has the same rational homology as $S^2$, and $H_2(C_2(R);\QQ)$ has a canonical generator $[S]$
that is the homology class of a fiber of $U\check{R} \subset C_2(R)$, oriented as the boundary of the unit ball of a fiber of $T\check{R}$.
For a $2$-component link $(J,K)$ of $\check{R}$, the homology class $[J \times K]$ of the image of
$J \times K$ in $H_2(C_2(R);\QQ)$ reads $lk(J,K)[S]$, where $lk(J,K)$ is the linking number of $J$ and $K$, see \cite[Proposition 1.6]{lesmek}.

\subsection{On propagators}
\label{subprop}

When $R$ is a rational homology sphere, a {\em propagator\/} of $C_2(R)$ is a $4$--cycle $F$ of $(C_2(R),\partial C_2(R))$ that is Poincar\'e dual to the preferred generator of
$H^2(C_2(R);\QQ)$ that maps $[S]$ to $1$.
For such a propagator $F$, for any $2$-cycle $G$ of $C_2(R)$,
$$[G]=\langle F,G\rangle_{C_2(R)} [S]$$
in $H_2(C_2(R);\QQ)$.

Let $B$ and $\frac12 B$ be two balls in $\RR^3$ of respective radii $\ell$ and $\frac{\ell}2$, centered at the origin in $\RR^3$.
Identify a neighborhood of $\infty$ in $R$ with $S^3 \setminus (\frac12 B)$ in $(S^3=\RR^3 \cup \{ \infty\})$
so that $\check{R}$ reads $\check{R}=M \cup_{]\ell/2,\ell] \times S^2} (\RR^3 \setminus (\frac12 B))$ for a rational homology ball $M$ whose complement in $\check{R}$ is identified with $\RR^3 \setminus B$.
There is a canonical regular map
$$p_{\infty} \colon (\partial C_2(R) \setminus UM) \rightarrow S^2$$
that maps the limit in $\partial C_2(R)$ of a sequence of ordered pairs of distinct points of $(\check{R} \setminus M)$ to the limit of the direction from the first point to the second one. See \cite[Lemma 1.1]{lesconst}.
Recall that $\tau_s\colon \RR^3 \times \RR^3 \rightarrow T\RR^3$ denotes the standard parallelization of $\RR^3$. Also recall that the sections $X$ of $UM$ that we consider are {\em constant\/} on $\partial M $, i.e. they read $\tau_s(\partial M \times \{V(X)\})$ for some fixed $V(X) \in S^2$ on $\partial M$. Let $X$ be such a section.
Then the {\em propagator boundary\/} $b_{X}$ associated with $X$ is the following $3$--cycle of $\partial C_2(R)$
$$b_{X}=p_{\infty}^{-1}(V(X)) \cup X(M)$$
and a {\em propagator associated with the section $X$\/} is a $4$--chain $F_{X}$ of $C_2(R)$ whose boundary reads
$b_{X}$. Such an $F_{X}$ is indeed a propagator because the algebraic intersection in $UM$ of a fiber and the section $X(M)$ is one.

\subsection{On the $\Theta$-invariant of a combed $\QQ$-sphere}
\label{subTheta}

\begin{theorem}\label{thmdefinvcomb}
Let $X$ be a section of $UM$ (that is constant on $\partial M$) for a rational homology ball $M$, and let $(-X)$ be the opposite section.
Let $F_{X}$ and $F_{-X}$ be two associated transverse propagators. Then $F_{X} \cap F_{-X}$ is a two-dimensional cycle whose homology class is independent of the chosen propagators. It reads
$\Theta(M,X)[S]$, where $\Theta(M,X)$ is a rational valued topological invariant of $(M,[X])$.
\end{theorem}
\bp
Recall that $C_2(R)$ has the same rational homology as $S^2$.
In particular, since $H_3(C_2(R);\QQ)=0$, there exist transverse propagators $F_{X}$ and $F_{-X}$ with the given boundaries $b_{X}$ and $b_{-X}$.
Without loss, assume that $F_{\pm X} \cap \partial C_2(R)=b_{\pm X}$.
Since $b_{X}$ and $b_{-X}$ do not intersect, $F_{X} \cap F_{-X}$ is a $2$--cycle. Since $H_4(C_2(R);\QQ)=0$, the homology class of $F_{X} \cap F_{-X}$ in $H_2(C_2(R);\QQ)$ does not depend on the choices of $F_{X}$ and $F_{-X}$ with their given boundaries. It reads
$\Theta(M,X)[S]$. Then it is easy to see that $\Theta(M,X) \in \QQ$ is a locally constant function of the section $X$.
\eop

When $R$ is an integral homology sphere, a combing $X$ is the first vector of a unique parallelization $\tau(X)$ that coincides with $\tau_s$ outside $M$, up to homotopy. When $R$ is a rational homology sphere, and when
$X$ is the first vector of a such a parallelization $\tau(X)$, this parallelization is again unique.
In this case, according to \cite[Section 6.5]{lessumgen} (or \cite[Theorem 2.14]{lesbonn}), the invariant $\Theta(M,X)$ is the degree $1$ part of the Kontsevich invariant of $(M,\tau(X))$ \cite{ko,kt,lesconst} and 
$$\Theta(M,X)=6 \lambda(M) + \frac{p_1(\tau(X))}{4}.$$

With our extension of the definition of $p_1$ to combings, we prove that the above formula also holds for combings.

\begin{theorem}
\label{thmvarinvcomb}
Let $X$ and $Y$ be two transverse sections of $UM$.
Then $$\Theta(M,Y)-\Theta(M,X)= lk(L_{X=Y},L_{X=-Y}).$$
In particular, $$\Theta(M,X)=6 \lambda(M) + \frac{p_1(X)}{4}.$$
\end{theorem}
\bp Let us prove that $\Theta(M,Y)-\Theta(M,X)= lk(L_{X=Y},L_{X=-Y})$.
This can be done as follows. Let $F_{-1}(\pm X,\pm Y)$ be the chain $F(\pm X,\pm Y)$ of Lemma~\ref{lemFXY} translated by $-1$ and seen in a collar $[-1,0] \times UM$ of $UM$ in $C_2(R)$. Assume that $F_X$ and $F_{-X}$ behave as products $[-1,0] \times \partial F_{\pm X}$ in $[-1,0] \times UM$. Then replacing these parts by $F_{-1}(X,Y)$ and $F_{-1}(-X,-Y)$, respectively, and making the appropriate easy corrections in $C_2(R)\setminus C_2(M)$ transforms $F_X$ and $F_{-X}$ into chains $F_Y$ and $F_{-Y}$ so that
$[F_Y \cap F_{-Y}]=[F_X \cap F_{-X}] + [F_{-1}(X,Y)\cap F_{-1}(-X,-Y)] $ where $[F_{-1}(X,Y)\cap F_{-1}(-X,-Y)]=lk(L_{X=Y},L_{X=-Y})[S]$ according to Proposition~\ref{propFXYcap}.
\eop

\def\cprime{$'$}
\providecommand{\bysame}{\leavevmode ---\ }
\providecommand{\og}{``}
\providecommand{\fg}{''}
\providecommand{\smfandname}{\&}
\providecommand{\smfedsname}{\'eds.}
\providecommand{\smfedname}{\'ed.}
\providecommand{\smfmastersthesisname}{M\'emoire}
\providecommand{\smfphdthesisname}{Th\`ese}

\newpage
\section*{Changes with respect to the first version v1}
There was an incorrect ``homotopy'' in the last sentence of the proof of Lemma 2.4 in v1. I apologize.
It is replaced by ``modification'' in the proof of the corresponding Lemma 2.1 in v2.\\
References have been added, thanks to Patrick Massot's help.\\
The paragraph order has been modified as follows:\\
$$\begin{array}{lll}
\mbox{v1}.\mbox{End of }\S 1.2 &\sim & \mbox{v2}. \mbox{Definition 3.5}\\ 
\mbox{v1}.\S 2.1 & \sim & \mbox{v2}.\S 4.1\\
\mbox{v1}.\S 2.2-\S 2.4 & \sim & \mbox{v2}.\S 2.1-\S 2.3\\
\mbox{v1}.\mbox{(Proof of Thm 1.1 in \S 3.2)} & \sim & \mbox{v2}.\S 3.2 \\
\mbox{v1}. \mbox{(\S 3.2 except Proof of Thm 1.1 + Proof of Thm 1.3 in \S 4.1)} & \sim & \mbox{v2}.\S 4.2\\
\mbox{v1}.\S 4.1 &\mbox{into} &\mbox{v2}.\S 4.4\\
\mbox{v1}.\S 4.2 & = &\mbox{v2}.\S 4.3\\
\mbox{v1}.\mbox{Definition 2.6} &=& \mbox{v2}.\mbox{Definition 4.16}\\
\end{array}$$
The redaction of v2.\S 2.1 has also been modified, and there are minor other local changes in the redaction elsewhere.\\
\end{document}